\documentclass{article}
\usepackage[T1]{fontenc}
\usepackage[utf8]{inputenc}
\usepackage{morefloats}
\usepackage[letterpaper]{geometry}
\usepackage{float}
\usepackage{amsmath}

\usepackage{graphicx}
\usepackage{color}
\usepackage{epstopdf}
\usepackage{adjustbox}
\usepackage{diagbox}
\usepackage{multirow}
\usepackage{array}
\usepackage{mathrsfs}
\usepackage{amsmath, amssymb, amsthm}
\usepackage{subfigure}
\usepackage{subfiles}
\usepackage{tikz}

\graphicspath{{./Figures2/}}

\newtheorem{theorem}{Theorem} 

\newtheorem{lemma}[theorem]{Lemma}

\usepackage{amsfonts} 
\usepackage{color}
\usepackage{adjustbox}
\usepackage{diagbox}
\definecolor{black}{rgb}{0,0,0}

\definecolor{red}{rgb}{1,0,0}

\usepackage{multirow}

\makeatother
\newcommand{\rR}{\mathrm{R}}
\newcommand{\bx}{\mathbf{x}}
\newcommand{\bu}{\mathbf{u}}
\newcommand{\bv}{\mathbf{v}}
\newcommand{\bw}{\mathbf{w}}
\newcommand{\bz}{\mathbf{z}}
\newcommand{\bn}{\mathbf{n}}
\newcommand{\cA}{\mathcal{A}}
\newcommand{\cT}{\mathcal{T}}
\newcommand{\cE}{\mathcal{E}}
\newcommand{\tilV}{\tilde{V}}
\usepackage{babel}
\usepackage{authblk}
\title{}

\title{\textbf{A conservative multiscale method for stochastic highly heterogeneous flow}}
 \author{Yiran Wang\thanks{Department of Mathematics, The Chinese University of Hong Kong, Shatin, Hong Kong SAR.}, \;
 	Eric Chung\thanks{Department of Mathematics, The Chinese University of Hong Kong, Shatin, Hong Kong SAR.} \;
 	and \; Shubin Fu\thanks{Department of Mathematics, University of Wisconsin-Madison, WI, USA. Corresponding author (shubinfu89@gmail.com)}
 }

\begin{document}
\maketitle

\begin{abstract}
In this paper, we propose a local model reduction approach for  subsurface flow problems in stochastic and highly heterogeneous media. To guarantee the mass conservation, we consider the mixed formulation of the flow problem and aim to solve the problem in a coarse grid to reduce the complexity of a large-scale system.
We decompose the entire problem into a training and a testing stage, namely the offline coarse-grid multiscale basis generation stage and online simulation stage with different parameters.
In the training stage, a parameter-independent and small-dimensional multiscale basis function space is constructed, which includes the media, source and boundary information. The key part of the basis generation stage is to solve some local problems defined specially.
With the parameter-independent basis space,  one can efficiently solve the concerned problems corresponding to different samples of permeability field in a coarse grid without repeatedly constructing a multiscale space for each new sample. A rigorous analysis on  convergence of the proposed method is proposed. In particular, we consider a generalization error, where bases constructed with one source will be used to a different source. In the numerical experiments, we apply the proposed method for both single-phase and two-phase flow problems.
Simulation results for both  2D and 3D representative models demonstrate the 
high accuracy and impressive performance of the proposed model reduction techniques.
\end{abstract}

\section{Introduction}
Many problems of fundamental and practical importance have multiple scales and high contrast properties. For instance, porous media, sedimenting suspensions and fluidized beds exhibit multiscale nature. Solving these problems in a fine grid that is sufficient to capture all small scale information is numerically expensive. Besides, when uncertainty is incorporated in physical models as well as numerical simulations, one needs to parameterize the partial differential equations with some random variables. Using sampling methods, one needs to call a deterministic solver once for each sample. Hence, developing an efficient deterministic fast solving algorithm is indispensable especially when we need to solve the problems repeatedly with a large number of samples. Because of the above considerations, some model reduction methods are necessary. In the past several decades, many researchers proposed various model reduction techniques, for example, upscaling methods \cite{barker1997critical,iliev2010fast,durlofsky1998coarse} and homogenization \cite{Chechkin.2007}. These methods aim to obtain a reduced model and one can therefore solve the problem in a coarse grid. However, some important small scale information is hard to be incorporated in the reduced models, which may hinder producing an accurate approximation. Multiscale methods \cite{MHM2,online-dg,Efendiev_GKiL_12,Aarnes2005257,ceg10,vasilyeva2019multiscale,chung2020generalized}, on the other hand,
are motivated to capture small-scale information of the underlying heterogeneous media in multiscale bases. Even using only a limited number of multiscale bases, the reduced-order solution  can still have relatively good accuracy. In other words, multiscale methods do well in attaining a trade-off between accuracy and computation costs.
Although one need first generate multiscale basis functions before the online simulation, one could efficiently solve the problem on a coarse grid and thus save online computational cost substantially. 
In addition to the spatial heterogeneities, here we also aim to construct coarse basis functions that are independent of the stochastic parameters so we do not need to recompute the multiscale bases. 

In  subsurface modeling, local mass conservation is vitally important for the transportation of the solute. To this end, mixed multiscale finite element methods \cite{aarnes04,chen2003mixed,aarnes2008mixed,chung2015mixed}, mortar multiscale methods \cite{ mortaroffline, ArPeWY07,Aarnes2005257,arbogast2006subgrid}, finite volume methods \cite{Jiang_Mish_MSFV_12,wang2016monotone,hajibeygi2008iterative,xu2011point} and some kinds of post-processing approaches \cite{odsaeter2017postprocessing,bush2013application,wang2021online} are proposed.
Among these mass conservative multiscale methods, the mixed multiscale finite element method (MMsFEM) \cite{chen2003mixed} has been successfully applied for various types of flow simulations. In MMsFEM, one velocity multiscale basis is constructed via solving a locally defined flow problem while  piecewise constant 
functions defined in each coarse-element are utilized for the pressure.
Although the MMsFEM succeeds in many scenarios, it fails to deal with very complicated porous media, and this motivates the development of the 
mixed  generalized multiscale finite element method (MGMsFEM) \cite{chung2015mixed,wang2021comparison}. In MGMsFEM, multiple multiscale basis functions for velocity are constructed and thus can capture more complicated media information. 

The construction of generalized multiscale space contains two steps. First, one needs to solve a series of local problems with a zero Neumann boundary condition to obtain some local snapshot spaces that include all possible local solutions. 
We then extract the dominant modes of the rich snapshot space by solving well-designed spectral problems which are motivated by analysis.
 The final multiscale basis function space is the direct sum of all  local multiscale spaces. It is worth mentioning that the construction of local basis function in different local regions are independent of each other, therefore parallel computation 
 can be easily adopted and thus these stages are cheap. These multiscale basis functions include important local permeability information  and can generate a good 
 coarse-grid solution without using too many computing resources. However, it does not include global media and source information and this brings about a precision limitation \cite{chung2015residual, online-mixed, online-dg}. 
The residual-driven multiscale basis functions were then proposed 
\cite{chung2015residual, online-mixed, online-dg} to improve the  accuracy. These residual-driven multiscale basis functions 
are computed via solving local problems using local residuals which means they 
include global media information, boundary condition information and 
source information. Besides, one can iteratively construct multiple local 
residual-driven bases and it is shown that with only a small number of these bases  one can obtain highly accurate coarse-grid solutions. 
The oversampling technique \cite{yang2020online}  can be utilized to further boost the performance of coarse-grid approximation. We will follow the main 
steps in \cite{yang2020online} to 
generate efficient multiscale basis functions for coarse-grid simulations of the stochastic flow problems. 
The basis space will include both the local permeability dependent multiscale 
basis functions and residual-driven multiscale basis functions. Both these two types of multiscale basis functions are constructed before the repeated online simulations and can  be parallelized without too many difficulties.

To deal with the uncertainty in the media, we use the Karhunen-Lo$\grave{e}$ve (KL) expansion \cite{huang2001convergence} in a reversed fashion. In particular, one can parameterize the stochastic permeability fields with a set of random variables. This process is divided into three steps. First, we choose a covariance function. Then, a spectral problem is solved to obtain some eigenfunctions corresponding to the dominating eigenvalues. The final stochastic field is generated from combinations of the selected eigenfunctions, where the same number of random variables serve as coefficients. Since a minority of eigenfunctions can contain most energy, one can obtain a low-rank representation of stochastic field. Based on this design, a sample permeability field will be straightly attained using a sample coefficient. Even though the number of samples is significantly reduced, repeatedly construction of multiscale bases for each sample is not judicious. To handle this challenge, many efforts 
can be found in literature \cite{jiang2016reduced,WANG2022114688,jin2011eulerian,hou2017exploring,zhang2015multiscale,hou2019model,li2020data,chen2021low,chen2020randomized,lei2015constructing,bright2016classification,ou2020low}. Here we apply the MGMsFEM with a representative permeability field to construct a parameter independent multiscale space. The involved field is also called the training permeability field. Once the multiscale space is constructed, one could make use of the approximation space to cheaply solve the flow problems corresponding to different samples on a coarse grid. 
 We rigorously analyze the proposed method by deriving a bound for the generalization error. We show that this error is a consequence of different sources and permeability fields between training and testing stages. To verify theoretical convergence results, we test the proposed method for both single-phase flow  and two-phase flow simulations. 
For the two-phase model, we use the finite volume method to solve the cheap transport equation on a fine grid and 
 efficiently solve the elliptic problem on a coarse grid with parameter independent 
 multiscale basis functions.
 We will predict the water cut in the two-phase model, which is a significant quantity in real applications. We demonstrate the efficiency of our methods using five benchmark models, where both 2D and 3D media are considered. 
The observed error decay behavior confirms our theoretical analysis and our method 
can fast generate reliable coarse-grid solutions.

The paper is organized as follows. We first present some preliminaries in Section 2 and then describe the construction of multiscale basis functions in details in Section 3.  Section 4 is devoted to analyze the proposed local model reduction method. The numerical experiments are provided in Section 5, where we will present simulation results for both the single-phase and two-phase flow models.  The paper is concluded in last Section.

\section{Preliminaries}
In this section, we will give some notations that will be used in the following presentations.
We consider the following problem in a given domain $\Omega\in \rR^d$, $d=2,3$ and a sample space $\Omega_r$:
\begin{align}
	\kappa(\bx;\omega)^{-1}\bv+\nabla p&=0,\quad \bx\in \Omega,\quad \omega\in \Omega_r,\label{model_v}\\
	\text{div}(\bv)&=f(\bx), \quad \bx\in \Omega,\label{model_p}
\end{align}
with zero Neumann boundary condition $\bv\cdot \bn=0.$ Here $\kappa(\bx;\omega)$ is a possibly highly heterogeneous random permeability field with high contrast. In particular, define $\kappa_{\text{contra}}=\frac{\max_{\bx\in \Omega}\kappa(\bx;\omega)}{\min_{\bx\in \Omega}\kappa(\bx;\omega)}$ and $\kappa_{\text{contra}}\gg 1$. We let $f$ be a given source term.

\begin{figure}[!htbp]
	\centering
	\includegraphics[width=0.9\textwidth]{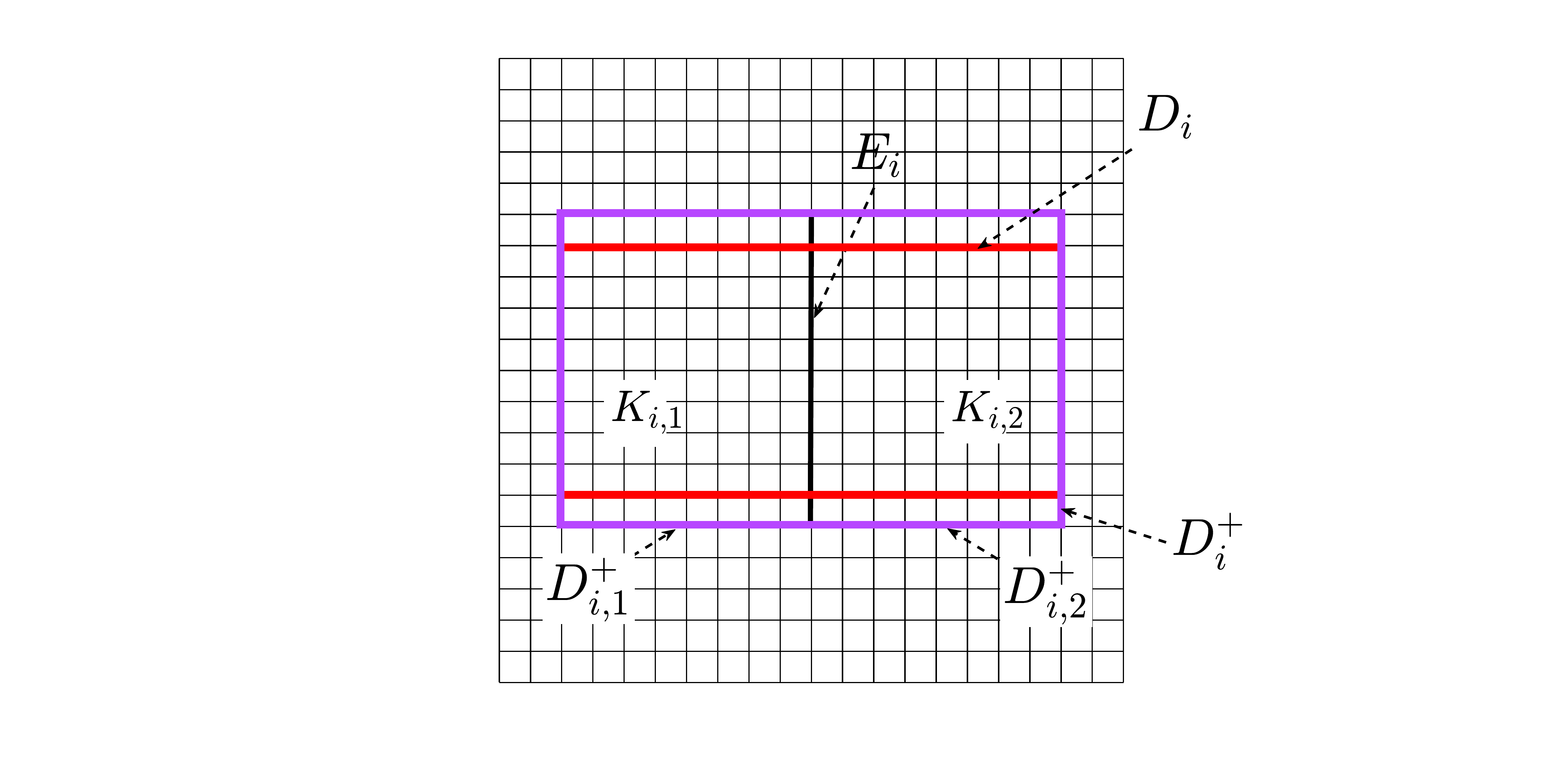}
	\caption{Illustration of a coarse-scale edge $E_i$, a corresponding coarse neighborhood $D_i$, and an oversampling neighborhood $D_i^+$.}
	\label{fig:mesh}
\end{figure}

First of all, we define the following inner products and norms.
\begin{align}
	\langle p,q\rangle=\int_{\Omega} pq, &\quad \cA(\bu,\bv)=\int_{\Omega} (\kappa(\bx;\omega))^{-1}\bu\cdot \bv,\\
	\|p\|_{L^2}^2=\int_{\Omega} |p|^2,&\quad \|\bu\|_a^2=\cA(\bu,\bu).
\end{align}
Then, we introduce the coarse and fine grids. Let $\cT_{H}$ be the coarse mesh which partitions $\Omega$ into a set of disjoint coarse-scale  elements $K_i$ with diameter $H$. In particular, we have $\bar{\Omega}=\cup_{i=1}^{N_{\text{e,c}}} K_i$, where $N_{\text{e,c}}$ is the number of coarse-scale  elements. Define $E_i$ to be a coarse-scale edge of coarse-scale element $K_i$ if $E_i=\partial K_i\cap K_j$ or  $E_c=\partial K_i\cap \partial\Omega$. Let $\cE_c(K_i)$ be the set of all the coarse-scale edges of $K_i$. Then we define $\cE_{c}=\cup_{i=1}^{N_{\text{e,c}}}\cE_{c}(K_i)$. 
We further define a coarse neighborhood $D_i$ to be a combination of two coarse-scale  elements sharing a common coarse-scale edge. More specifically, we define
\begin{align*}
	D_i:=\{K\in \cT_H: E_i\subset \partial K\}, \quad i=1,\ldots,N_{\text{E,c}},
\end{align*}
where $N_{\text{E,c}}$ is the number of coarse-scale edges. Moreover, we define an ovesampled coarse neighborhood corresponding to a coarse-scale edge $E_i$ by $D_i^+$, which is constructed by enlarging $D_i$ in some way. An example is shown in Figure \ref{fig:mesh}.  Based on the coarse grid $\cT_{H}$, we further construct a fine grid $\cT_h$. In particular, for each coarse-scale element $K_i$, we partition it into some fine-scale  elements with a fine mesh size $h$ and we denote the mesh by $\cT_h(K_i)$. Then $\cT_h=\cup_{i=1}^{N_{\text{e,c}}}T_h(K_i)$. Correspondingly, we define $\cE_f(K_i)$ to be all fine-scale edges in the partition $\cT_h(K_i)$. In addition, let $\cE_f^0(K_i)$ be the interior edges in $\cT_h(K_i)$. And $\cE_f=\cup_{i=1}^{N_{\text{e,c}}}\cE_f(K_i)$. We illustrate the above notations in Figure \ref{fig:mesh}.

The variational formulation of the exact problem \eqref{model_v}-\eqref{model_p} is: find $\bv\in (L^2(\Omega))^2$ and $p\in H_0^1(\Omega)$ such that 
\begin{align}
	\int_{\Omega}\kappa^{-1}\bv\cdot \bu-\int_{\Omega}\text{div}(\bu)p&=0,\quad \forall \bu\in (L^2(\Omega))^2,\\
	\int_{\Omega}\text{div}(\bv)q&=\int_{\Omega}fq,\quad \forall q\in H_0^1(\Omega),
\end{align}
with zero Neumann boundary condition $\bv\cdot \bn=0$. To solve $\bv$ and $p$ in a finite-dimensional space, we first introduce the fine-grid space. Let $V_f$ be the lowest Raviart Thomas vector field ($RT_0$) and $Q_f$ be space of piecewise constant functions on the fine mesh. We set $V_f=\text{span}\{\psi_{f,1},\cdots,\psi_{f,N_{\text{E,f}}}\}$ and
$Q_f=\text{span}\{p_{f,1},\cdots,p_{f,N_{\text{e,f}}}\}$, where $N_{\text{E,f}}$ and $N_{\text{e,f}}$ are the numbers of inner edges and blocks on the fine grid respectively. $V_{f}^0:=\{\bv\in V_{f}:\bv\cdot \bn=0 \text{ on }\partial \Omega\}$.
We obtain the fine-grid solution (reference solution) $(\bv_f,p_f)$ by solving the following system.
\begin{eqnarray}
	\begin{aligned}
		\cA(\bv_f,\bw)-\langle \text{div}(\bw),p_f\rangle&=0,\quad \forall \bw\in V_f^0,\\
		\langle \text{div}(\bv_f),q\rangle&=\langle f,q\rangle, \quad \forall q\in Q_f, \label{fine}
	\end{aligned}
\end{eqnarray}
where $\bv_f\cdot \bn=0$ on $\partial \Omega$. 

We further introduce the multiscale space $(V_{\text{ms}},Q_{\text{ms}})$. Here $Q_{\text{ms}}$ is the pressure multiscale space spanned by piecewise constant functions on each coarse-scale element, and $V_{\text{ms}}$ is the velocity multiscale space. To construct $V_{\text{ms}}$, we first construct a snapshot space $V_{\text{snap}}$ and then reduce it to a lower-dimensional approximation space $V_{\text{ms,0}}$ by some well-designed local spectral problems. Finally, we perform the residual-driven enrichment iterations on $V_{\text{ms,0}}$ to get the final multiscale space $V_{\text{ms}}$. We define $V_{\text{ms,k}}$ to be the multiscale space on the enrichment level $k$.

After we obtain the $(V_{\text{ms}},Q_{\text{ms}})$, we seek the multiscale solution $(\bv_{\text{ms}},p_{\text{ms}})$ by solving the following system.
\begin{eqnarray}
	\begin{aligned}
		\cA(\bv_{\text{ms}},\bu)-\langle\text{div}(\bu),p_{\text{ms}}\rangle&=0, \quad \forall \bu\in V_{\text{ms}},\\
		\langle\text{div}(\bv_{\text{ms}}),q\rangle&=\langle f,q\rangle,\quad \forall q \in Q_{\text{ms}},\label{coarse}
	\end{aligned}
\end{eqnarray}
where $\bv_{\text{ms}}\cdot \bn=0$. 
In the following part, we will rewrite \eqref{fine} and \eqref{coarse} in a matrix formulation. Let $v_{f,v}$ and $p_{f,v}$ to be the corresponding coefficient vector for reference solution $\bv_f$ and $p_f$. $A_f[i,j]=\cA(\psi_{f,i},\psi_{f,j})$ and $B_h[i,j]=\langle \text{div}(\psi_{f,i}),p_{f,j}\rangle$. Then we obtain $v_{f,v}$ and $p_{f,v}$ by solving
\begin{equation*}
	\left[\begin{array}{cc}
		A_{f}& B_{f}^T \\
		B_{f}  & 0
	\end{array}\right]\left[\begin{array}{l}
		v_{f,v} \\
		p_{f,v}
	\end{array}\right]=\left[\begin{array}{c}
		0 \\
		F_f
	\end{array}\right].
\end{equation*}
Similarly, let $v_{c,v}$ and $p_{c,v}$ to be the corresponding coefficient vectors for multiscale solution $\bv_{\text{ms}}$ and $p_{\text{ms}}$. Define $U$ to be the matrix that collects all the velocity multiscale solutions and $M_c$ to be the matrix that provides an embedding operation from $Q_{\text{ms}}$ to $Q_f$. $v_{c,v}$ and $p_{c,v}$ are obtained by solving
\begin{equation*}
	\left[\begin{array}{cc}
		U^TA_{f}U& U^T B_{f}^TM_c \\
		M_c^TB_{f}U  & 0
	\end{array}\right]\left[\begin{array}{l}
		v_{c,v} \\
		p_{c,v}
	\end{array}\right]=\left[\begin{array}{c}
		0 \\
		M_c^TF_f
	\end{array}\right].
\end{equation*}
In the following section, we will give a detailed review of the MGMsFEM.

\section{Mixed Generalized Multiscale Finite Element Method}
\begin{figure}[!htbp]
	\centering
	\includegraphics[width=0.9\textwidth]{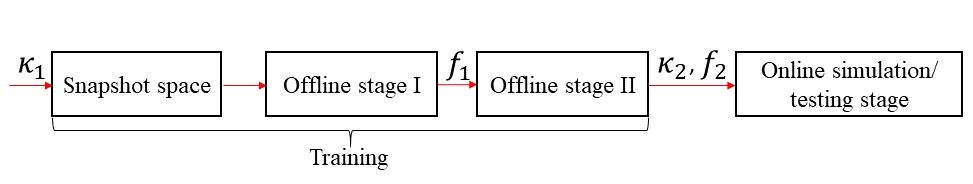}
	\caption{A training-testing algorithm. $\kappa_1$ and $\kappa_2$, $f_1$ and $f_2$ are permeability fields and sources used in training and testing stages. Training stage contains three steps and it is performed offline. The online simulation/ testing stage is conducted with different parameters $\kappa_2$ and $f_2$, which may be different with those used in the training stage.}\label{fig:algo}
\end{figure}

In this section, we will present the Generalized Multiscale Finite Element Method in mixed formulation (MGMsFEM).
As mentioned before, the construction of multiscale space is divided into three steps, which will be described in the following subsections. In particular, we first construct a set of local snapshot space. Then a reduced multiscale space is further constructed by solving local snapshot problems. The last step is to construct some residual-driven bases to enrich the multiscale space obtained in the second step. The last two steps are called offline  stage I and II. The corresponding multiscale bases in the above two stages are called multiscale basis type I and multiscale basis type II, where the latter one are residual-driven bases. The whole process of our algorithm is shown in Figure \ref{fig:algo}. Throughout this section, we take a deterministic permeability field $\kappa=\kappa(\bx,\bar{\xi})$, where $\bar{\xi}$ is a sample from the random space $\Omega_r$.
\subsection{Snapshot space}
In this subsection, a snapshot space will be constructed to represents all possible local solutions. First,  a set of local snapshot spaces $\{V_{\text{snap}}^{i}\}$ are spanned by local snapshot functions, which are solved in the fine grid. Then, the final snapshot space is defined by $V_{\text{snap}}=\bigoplus_{E_i\in \cE_{c}}V_{\text{snap}}^{i}$. 

Suppose $E_i\in \cE_c$ and $E_i$ can be partitioned into some fine-scale edges, i.e. $E_i=\cup_{j=1}^{L_i} e_j$. Here $L_i$ is the number of fine-scale edges on the coarse-scale edge $E_i$. To obtain $V_{\text{snap}}^i$, we solve the following system:
\begin{eqnarray}
	\begin{aligned}
		\kappa^{-1}\bv_j^i+\nabla p_j^i=0 \text{ in } D_i,\\
		\text{div}(\bv_j^i)=\beta_j^i \text{ in }D_i. \label{local snap}
	\end{aligned}
\end{eqnarray}
subject to the zero Neumann boundary condition $\bv_j^i\cdot n_i=0$ on $\partial D_i$, where $n_i$ is a fixed unit-normal vector for $E_i$. We solve \eqref{local snap} separately in $K_{i,1}$ and $K_{i,2}$ (see the Figure \ref{fig:mesh} for illustration). Hence we use an additional boundary condition $\bv_j^i\cdot n_i=\delta_j^i$, where $\delta_j^i$ is the Kronecker delta. In particular, it is defined as follows.
\begin{equation*}
	\delta_j^i=\left\{
	\begin{array}{cc}
		1,&\text{ on } e_j,\\
		0,& \text{ on }E_i\setminus e_j,
	\end{array}
	j=1,\ldots,L_i.
	\right.
\end{equation*}
Here, a compatible condition $\int_{K_{i,p}} \beta_j^i=\int_{\partial K_{i,p}} \bv_j^i\cdot n_i$, $p=1,2$. We solve the above local problem in each $D_i$ to generate $V_{\text{snap}}^i$. Finally, we form the snapshot space by using $V_{\text{snap}}=\bigoplus_{E_i\in \cE_c} V_{\text{snap}}^i$.
\subsection{Offline stage I}
In this subsection, we perform a dimension reduction on $V_{\text{snap}}$, whose dimension is comparable to  that of fine-grid space $V_f$. The resulted small-dimensional space is referred as the multiscale space in the offline  stage I. 
First of all, we define two symmetric and positive definite bilinear operators $a(\cdot,\cdot)$ and $s(\cdot,\cdot)$ on $V_{\text{snap}}^i\times V_{\text{snap}}^i$. More specifically, 
\begin{eqnarray}
	\begin{aligned}
		a(\bv,\bu)&=\int_{E_i}\kappa^{-1}(\bv\cdot \bn_i)(\bu\cdot \bn_i),\\
		s(\bv,\bu)&=\frac{1}{H}\left(\int_{D_i}\kappa^{-1}\bv\cdot \bu+\int_{D_i}\text{div}(\bv)\text{div}(\bu)\right),
	\end{aligned}
\end{eqnarray}
for $\bv,\bu\in V_{\text{snap}}^i$, and $\bn_i$ is the fixed unit normal vector for $E_i$.
We solve the following local spectral problem to get $\lambda\in \rR$ and $\bv\in V_{\text{snap}}^i$,
\begin{equation}
	a(\bv,\bu)=\lambda s(\bv,\bu),\quad \forall \bu\in V_{\text{snap}}^i.\label{spectral}
\end{equation}
After solving \eqref{spectral} in $D_i$, we arrange the eigenvalues in ascending order $\lambda_{1}^i\leq \lambda_2^i\leq \cdots \leq \lambda_{L_i}^i$. Let $\psi_j^i$ be the eigenfunction corresponding to $\lambda_j^i$. To reduce the dimension, we here keep the eigenfunctions associated with $l_i$ smallest eigenvalues and use them to span the local offline space $V_{\text{ms}}^{i}$, i.e. $V_{\text{ms}}^{i}=\text{span}\{\psi_1^i,\ldots,\psi_{l_i}^i\}$. Then the global offline space is $V_{\text{ms}}=\bigoplus_{E_i\in \cE_c}V_{\text{ms}}^{i}$.

The eigenvalue problem \eqref{spectral} can be recast in the following matrix formulation.
Define $A_{\text{snap}}^{i}$ and $S_{\text{snap}}^{i}$ as follows.
\begin{align}
	A_{\text{snap}}^{i}[k,j]=a(\bv_k^i,\bv_j^i),\quad S_{\text{snap}}^{i}[k,j]=s(\bv_k^i,\bv_j^i), \quad \forall \bv_k^i,\bv_j^i\in V_{\text{snap}}^i.
\end{align}
Define $\bv_{f,v}$ to be the coefficient vector of the solution $\bv$ to \eqref{spectral} w.r.t the snapshot basis functions in $V_{\text{snap}}^{i}$.
Then we can get $\bv_{f,v}$ and $\lambda$ by solving the following matrix formulation:
\begin{align*}
	A_{\text{snap}}^i\bv_{f,v}=\lambda S_{\text{snap}}^{i}\bv_{f,v}.
\end{align*}
\subsection{Offline stage II}
The multiscale basis function constructed in last section only depends on local permeability field{\color{red}}, where no source or boundary condition of the solution is included in the basis construction. To compensate this deficiency, one can construct new bases to enhance the approximation ability of multiscale space, which is the goal of this subsection. We recall that the key idea in the offline stage I is utilizing a set of  spectral problems to reduce the dimension of approximation space. However, since our motivation here is to perform some enrichments effectively, we apply a different bases construction method.  More specifically, we will enrich the multiscale space from offline stage I by some residual-driven bases. To this end, we first solve \eqref{local snap} in some oversampled regions as follows.  Figure \ref{fig:mesh} shows an oversampled domain. For each $i=1,\ldots, N_{E,c}$, we solve
\begin{eqnarray}
	\begin{aligned}
		\kappa^{-1}\bv_j^{i,+}+\nabla p_j^{i,+}=0 \text{ in } D_i^+,\\
		\text{div}(\bv_j^{i,+})=\gamma_j^i,\text{ in }D_i^+ \label{oversample snap}
	\end{aligned}
\end{eqnarray}
subject to the zero boundary condition $\bv_j^{i,+}\cdot n_i=0$ on $\partial D_i^{+}$. Similarly, we separately solve \eqref{oversample snap} in $D_{i,1}^+$ and $D_{i,2}^+$ (Figure \ref{fig:mesh}) with the additional boundary condition $\bv_{j}^{i,+}\cdot \bn_i=\delta_{j}^i$ on $E_i^+$. Here $E_i^+$ is an oversampled coarse-scale edge shared by $D_{i,1}^+$ and $D_{i,2}^+$. Moreover, the constant $\gamma_j^i$ is selected such that \eqref{oversample snap} is solvable. From \eqref{oversample snap}, we can obtain $\bv_{1}^{i,+},\ldots,\bv_{L_i^+}^{i,+}$ to span the oversampled local snapshot space $V_{\text{snap}}^{i,+}$. Let $\tilV_{\text{snap}}^{i,+}$ be the divergence free subspace of $V_{\text{snap}}^{i,+}$. In the following part, we seek the residual-driven bases in the oversampled snapshot space. 

Suppose $D\subset \Omega$ is a specific domain. Let $V_{D}$ be the subspace of $V_{\text{snap}}$ that contains functions supported in $D$. In particular, we have $V_{D}=\bigoplus_{D_i\subset D}V_{\text{snap}}^{i}$. We further define $\tilde{V}_{D}$ to be the divergence free subspace of $V_D$. Define the residual operator $R_{D}$ on $V_D$ by
\begin{align*}
	R_{D}(\bv)=\int_{D}\kappa^{-1}\bv_{\text{ms}}\cdot \bv-\int_{D}\text{div}(\bv)p_{\text{ms}},\quad \forall \bv\in V_{D}.
\end{align*}
If we restrict $R_{D}$ on $\tilde{V}_{D}$, we have
\begin{align*}
	R_{D}(\bv)=\int_{D}\kappa^{-1}\bv_{\text{ms}}\cdot \bv,\quad \forall \bv\in \tilde{V}_{D}.
\end{align*}
We denote the operator norm of the residual $R_{D}$ by $\|R_{D}\|$. The construction of residual-driven bases is performed by iterations. We denote the multiscale space after $n$ iterations by $V_{\text{ms},n}$. In particular, $V_{\text{ms},0}$ is the offline multiscale space. Let $k$ be the iteration level. We start at $k=0$. Each iteration contains five steps as follows. 

Step 1: Compute the multiscale solution $(\bv_{\text{ms}}^k,p_{\text{ms}}^k)$ in the current multiscale space $V_{\text{ms,k}}\times Q_{\text{ms}}$. In particular, we solve
\begin{eqnarray}
	\begin{aligned}
		\int_{\Omega}\kappa^{-1}\bv_{\text{ms}}^k\cdot \bu-\int_{\Omega}\text{div}(\bu)p_{\text{ms}}^k&=0,\quad \forall \bu\in V_{\text{ms,k}},\\
		\int_{\Omega}\text{div}(\bv_{\text{ms}}^k)q&=\int_{\Omega}fq,\quad \forall q\in Q_{\text{ms}}.
		\label{ms_k}
	\end{aligned}
\end{eqnarray}
Suppose $V_{\text{ms,k}}=\{\psi_1,\ldots,\psi_{N_k}\}$, where $N_k$ is the number of multiscale bases in $V_{\text{ms,k}}$. $Q_{\text{ms}}=\{p_1,\ldots,p_{\text{e,c}}\}$.
We can rewrite \eqref{ms_k} in a matrix formulation as follows,
\begin{equation*}
	\left[\begin{array}{cc}
		A_{\text{ms,k}}& B_{\text{ms,k}}^T \\
		B_{\text{ms,k}}& 0
	\end{array}\right]\left[\begin{array}{l}
		\bv_{\text{ms,v}}^k \\
		p_{\text{ms,v}}^k
	\end{array}\right]=\left[\begin{array}{c}
		0 \\
		F_{\text{ms,k}}
	\end{array}\right].
\end{equation*}
Here, $A_{\text{ms,k}}$ and $B_{\text{ms,k}}$ are constructed using  multiscale bases in $V_{\text{ms},k}$. In particular, $A_{\text{ms,k}}[i,j]=\mathcal{A}(\psi_i,\psi_j)$ and $B_{\text{ms,k}}[i,j]=\langle\text{div}(\psi_i),p_j \rangle$.

Step 2: Determine whether a enrichment is needed by some particular error indicators. In particular, given an error tolerance $\tau$, if $\|R_{\Omega}\|>\tau$, we go through the following three steps.

Step 3: Select a set of local regions to perform the residual-driven enrichment. For each $E_i\in \cE_{c}$, we construct a corresponding oversampled neighborhood $D_i^+$. See the Figure \ref{fig:mesh} for an example of an oversampled neighborhood. Let $D_{1},\ldots,D_{p_k}$ be a set of local neighborhoods that form a non-overlapping partition of $\Omega$.

Step 4: Obtain the boundary condition for the residual-driven bases. For each $D_i^+$, we solve for $\psi_i^+\in \tilV_{\text{snap}}^{i,+}$ such that
\begin{align}
	\int_{D_i^+}\kappa^{-1}\psi_i^+\cdot \bu=R_{D_i^+}(\bu), \quad \forall \bu\in \tilV_{\text{snap}}^{i,+}.
\end{align} 
Since $\psi_i^+$ is supported in $D_i^+$, we need to restrict $\psi_{i}^+\cdot \bn_i^+$ on the coarse-scale edge $E_i$ and normalize it to obtain $z_i$.

Step 5: For each $D_i$, compute a residual-driven basis $\phi_i$ supported in $D_i$ with computed boundary value $z_i$ by solving
\begin{eqnarray}
	\begin{aligned}
		\kappa^{-1}\phi_i+\nabla p_i&=0\text{ in } D_i,\\
		\text{div}(\phi_i)&=\beta_i\text{ in } D_i,\\
		\phi_i\cdot \bn_i&=z_i\text{ on }E_i,\\
		\phi_i\cdot\bn_i&=0\text{ on } \partial D_i,\label{online}
	\end{aligned}
\end{eqnarray}
where $\beta_i$ is chosen to satisfy $\int_K\beta_i=\int_{\partial K}\phi_i\cdot \bn_i$ for each $K\subseteq D_i$ and $n_i$ is a fixed unit-normal vector for the coarse face $E_i$.
Note that the process of solving \eqref{online} is unnecessary when $V_{\text{snap}}$ is already obtained.
Define $\phi_i,\ldots,\phi_{p_k}$ to be residual-driven bases obtained in iteration $k$. We then update the velocity space by letting $V_{\text{ms}}^{k+1}=V_{\text{ms}}^k\cup \text{span}\{\phi_i,\ldots,\phi_{p_k}\}$.

Using the above algorithm, we can terminate the residual-driven enrichment if $\|R_{\Omega}\|\leq \tau$.

\begin{figure}[!htbp]
	\centering
	\includegraphics[width=1.\textwidth]{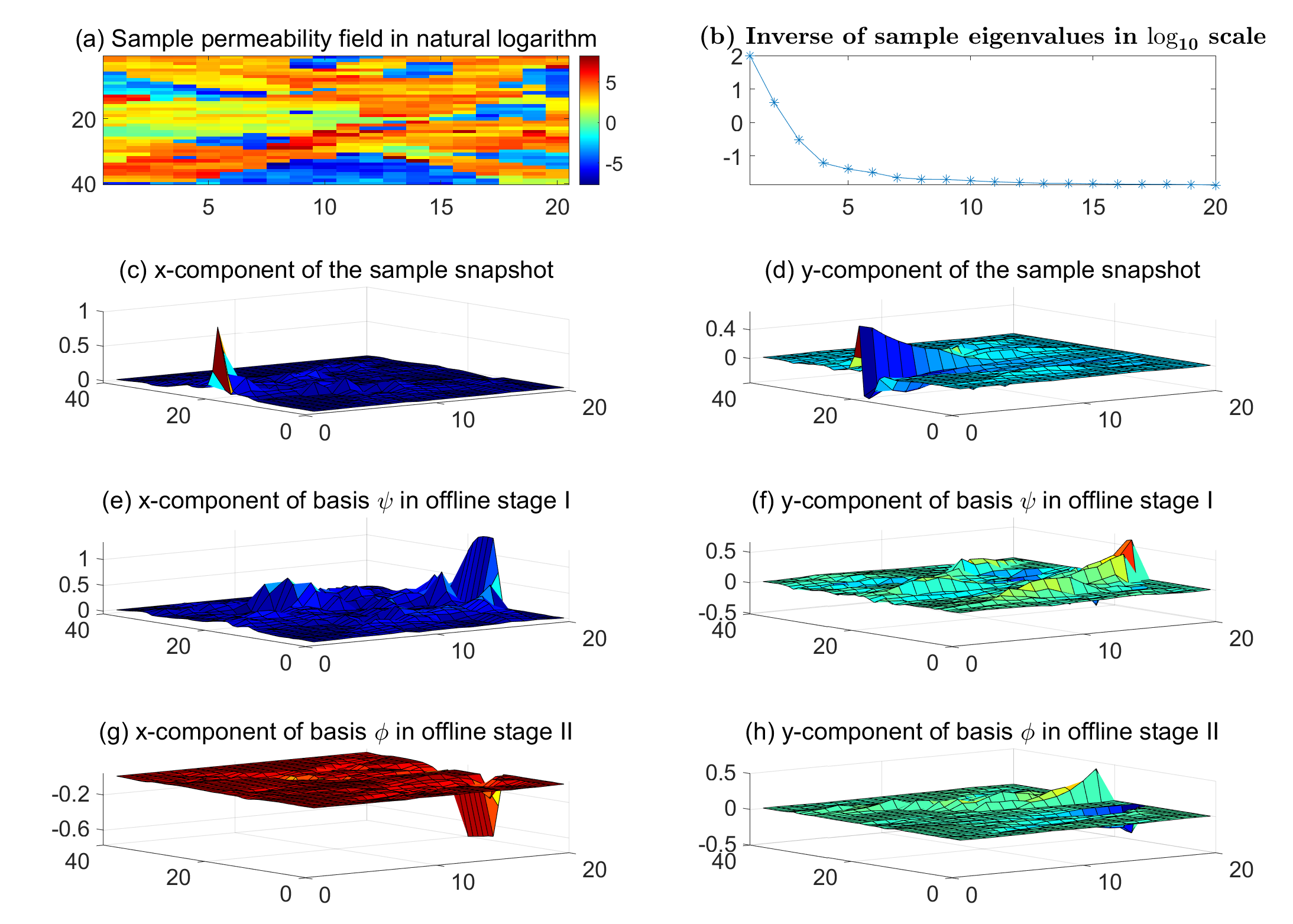}
	\caption{Sample local permeability field, inverse of sample eigenvalues, a sample snapshot, an offline basis $\psi$, and a residual-driven basis $\phi$ in a specific local neighborhood. (a): Sample permeability field in natural logarithm. (b): Inverse of 20 smallest eigenvalues shown in $\log_{10}$ scale. (c)-(d): $x$ and $y$ components of a sample local snapshot. (e)-(f): $x$ and $y$ components of a sample offline multiscale basis $\psi$ in offline  stage I. (g)-(h): $x$ and $y$ components of a sample residual-driven basis $\phi$ in the offline stage II.}
	\label{loc_fig}
\end{figure}

In Figure \ref{loc_fig}, local profile in a specific coarse neighborhood (see Figure \ref{fig:mesh}) is shown. Choosing a local region $D_i$, we plot the permeability field, a multiscale basis in the offline stage I, and a residual-driven basis corresponding to $D_i$. Moreover, a sample snapshot generated during the construction of multiscale bases is also displayed. From (a), we can see there are some high-contrast variations in the local permeability field. The highest value is over $3000$ while the smallest value is $10^{-4}$. The inverse of smallest eigenvalues are shown in $\log_{10}$ scale in (b). The decay is rapid especially among the 5 smallest eigenvalues. In (c) and (d), a sample snapshot is displayed, where $x$ and $y$ components are shown in two columns.  A sample multiscale basis $\psi$ in the offline  stage I and a sample residual-driven multiscale basis $\phi$ are displayed in (e)-(f) and (g)-(h) respectively. The multiscale bases show some distinct variations at some specific places due to the heterogeneity of the permeability field. Besides, the residual-driven basis function has fewer fluctuations compared with the type I basis function $\psi$ but the oscillations in $\psi$ are more centralized, where the distinction is a consequence of different constructions. In particular, the residual-driven basis function is computed based on some residuals while the basis type I is a linear combination of some eigenfunctions that are more divergent. 
We also remark that although the residual-driven basis is more efficient but it relies on 
the residual of a specific problem, so if the online simulation problem has a sharp difference with the training problem, residual-driven bases may not work well, and we will show this point in our analysis part. 

\section{Analysis}
Suppose $\Omega$ is a bounded domain in $R^d(d=2,3)$. $\Omega_r$ is a sample space. We aim at solving the proposed problem:
\begin{align}
	\kappa^{-1}(\bx;\omega)\bu+\nabla p&=0 \text{ in }\Omega ,\\
	\text{div}(\bu)&=f\text{ in }\Omega .\label{eq:model}
\end{align}
First of all, we introduce some notations. 
\begin{enumerate}
	\item $\bu_{f}(\bx;\omega,f)$: a fine-grid velocity solution solved in $V_f$. $p_{f}(\bx;\omega,f)$: a fine-grid pressure solution solved in $Q_f$.
	\item $\omega_1$ and $f_1$ are a training permeability field and a training source used to generate a multiscale space.
	\item $V_{\text{ms},k}$: a velocity multiscale space at enrichment iteration $k$. $V_{\text{ms}}:=V_{\text{ms},m}$: a final velocity multiscale space is generated after $m$ enrichment iterations.
	\item  $\bu_{\text{ms}}(\bx;\omega,f)$: a velocity multiscale solution solved in $V_{\text{ms}}$. $p_{\text{ms}}(\bx;\omega,f)$: a pressure multiscale  solution solved in $Q_{\text{ms}}$. $\bu_{\text{ms}}^{1,1,(k)}(\bx;\omega_1,f_1)$: a velocity multiscale solution corresponding to $\omega_1$ and $f_1$ solved in $V_{\text{ms},k}$.
	\item $\kappa_{1,\text{max}}:=\displaystyle\max_{\bx\in \Omega}\kappa(\bx,\omega_1)$,
	$\kappa_{1,\text{min}}:=\displaystyle\min_{\bx\in \Omega}\kappa(\bx,\omega_1)$,
	$\kappa_{1,\text{contra}}:=\frac{\kappa_{1,\text{max}}}{\kappa_{1,\text{min}}}$.
	\item $\Omega_1,\ldots,\Omega_P$ are a set of non-overlapping regions such that $\Omega=\bigcup_{j=1}^P \Omega_j$.
\end{enumerate}
In particular, $\bu_f(\bx;\omega,f)$ and $p_f(\bx;\omega,f)$ are velocity and pressure solutions to \eqref{eq:fine sol v} and \eqref{eq:fine sol p}; $\bu_{\text{ms}}(\bx;\omega,f)$ and $p_{\text{ms}}(\bx;\omega,f)$ are velocity and pressure solutions to \eqref{eq:ms sol v} and \eqref{eq:ms sol p} as follows.
\begin{align}
	\int_\Omega  \kappa^{-1}(\bx;\omega)\bu_f(\bx;\omega,f) \cdot\bv-\int_\Omega  \text{div}(\bv)p_f(\bx;\omega,f)&=0, \quad\forall \bv\in V_f^0, \label{eq:fine sol v}\\ 
	\int_\Omega  \text{div}(\bu_f(\bx;\omega,f))q&=\int_\Omega  fq, \quad\forall q\in Q_f.\label{eq:fine sol p}
\end{align}
\begin{align}
	\int_\Omega  \kappa^{-1}(\bx;\omega)\bu_{\text{ms}}(\bx;\omega,f)\cdot \bv-\int_\Omega  \text{div}(\bv)p_{\text{ms}}(\bx;\omega,f)&=0,\quad \forall \bv\in V_{\text{ms}},\label{eq:ms sol v}
	\\ 
	\int_\Omega  \text{div}(\bu_{\text{ms}}(\bx;\omega,f))q&=\int_\Omega  fq, \quad\forall q\in Q_{\text{ms}}.\label{eq:ms sol p}
\end{align}
In the following part, we give an  error estimate of a multiscale solution  in an online simulation stage solved in a fixed multiscale space $V_{\text{ms}}$. In particular, we arbitrarily choose a sample permeability field and a source term. Once they are chosen, they are deterministic, which are denoted by $\kappa_2$ and $f_2$ respectively. We recall that $\kappa_1$ and $f_1$ are the training permeability field and source used to generate $V_{\text{ms}}$.
To simplify the above mentioned notations, we use some abbreviations. For $i,j\in\{1,2\}$, we let $\bu_{f}^{i,j}:=\bu_{f}(\bx;\omega_i,f_j)$ and $p_{f}^{i,j}:=p_{f}(\bx;\omega_i,f_j)$; $\bu_{\text{ms}}^{i,j}:=\bu_{\text{ms}}(\bx;\omega_i,f_j)$ and $p_{\text{ms}}^{i,j}:=p_{\text{ms}}(\bx;\omega_i,f_j)$; $\bu_{\text{ms}}^{1,1,(m)}:=\bu_{\text{ms}}^{(m)}(\bx;\omega_1,f_1)$; $\kappa_i:=\kappa(\bx,\omega_i)$.

Define $e=\kappa_2^{-\frac{1}{2}}(\bu_{f}^{2,2}-\bu_{\text{ms}}^{2,2})$. To estimate $e$, we split it into five parts as follows:
\begin{eqnarray}
	\begin{aligned}
		e
		&=\kappa_2^{-\frac{1}{2}}\bu_{f}^{2,2}-\kappa_1^{-\frac{1}{2}}\bu_{f}^{1,2}+
		\kappa_1^{-\frac{1}{2}}(\bu_{f}^{1,2}-\bu_{f}^{1,1}+\bu_{f}^{1,1}-\bu_{\text{ms}}^{1,1}+\bu_{\text{ms}}^{1,1}-\bu_{\text{ms}}^{1,2})+\kappa_1^{-\frac{1}{2}}\bu_{\text{ms}}^{1,2}-\kappa_2^{-\frac{1}{2}}\bu_{\text{ms}}^{2,2}\\
		&:=e_1+e_2+e_3+e_4+e_5. \label{err_split}
	\end{aligned}
\end{eqnarray}
In particular,
\begin{align*}
	e_1&=\kappa_2^{-\frac{1}{2}}\bu_{f}^{2,2}-\kappa_1^{-\frac{1}{2}}\bu_{f}^{1,2},\\
	e_2&=\kappa_1^{-\frac{1}{2}}(\bu_{f}^{1,2}-\bu_{f}^{1,1}),\\
	e_3&=\kappa_1^{-\frac{1}{2}}(\bu_{f}^{1,1}-\bu_{\text{ms}}^{1,1}),\\
	e_4&=\kappa_1^{-\frac{1}{2}}(\bu_{\text{ms}}^{1,1}-\bu_{\text{ms}}^{1,2}),\\
	e_5&=\kappa_1^{-\frac{1}{2}}\bu_{\text{ms}}^{1,2}-\kappa_2^{-\frac{1}{2}}\bu_{\text{ms}}^{2,2}.
\end{align*}

We remark that our proposed method is to choose a well-chosen permeability field and a fixed source term, based on which a multiscale space is constructed in an offline stage. It is worth noting that the offline stage contains the process of computing residual-driven bases. Furthermore, we can use the pre-computed multiscale space to solve the problem with a different permeability field and a source term in a testing stage. Hence, the final error may come from three aspects: multiscale approximation, different permeability fields, and different source terms. More specifically, in \eqref{err_split}, $e_3$ refers to the approximation error of a multiscale space, i.e. we measure the difference of a multiscale solution with a corresponding reference solution, where the permeability field and the source term in the concerned equation are the same with those used in constructing multiscale space. $e_1$ and $e_2$ are associated with differences of reference solutions to different permeability fields and sources, respectively. $e_4$ and $e_5$ refer to differences of two mutiscale solutions. Besides, the way to derive the bounds for $e_1$ and $e_5$ is almost the same. Hence, it suffices to derive the error bounds for $e_1$.  At the same time, the bounds for $e_2$ and $e_4$ are obtained similarly. We will estimate each part individually and summarize them in the end. We will start with $e_1$ and $e_5$. First of all, we prove the stability of the fine-scale and multiscale velocity solutions.

\begin{lemma}\label{stab}
	Suppose $\bu_{f}:=\bu_{f}(\bx;\omega,f)$ and $\bu_{\text{ms}}:=\bu_{\text{ms}}(\bx;\omega,f)$ are fine-grid and multiscale velocity solutions to \eqref{eq:fine sol v}-\eqref{eq:fine sol p} and \eqref{eq:ms sol v}-\eqref{eq:ms sol p}, respectively. Then there exists $C$ (independent of $\kappa$) such that the following two estimations hold.
	\begin{align}
		\|\bu_{f}\|_{L^2}&\leq C\|f\|_{L^2},\label{stab_fine}\\
		\|\bu_{\text{ms}}\|_{L^2}&\leq C\|f\|_{L^2}.\label{stab_ms}
	\end{align}
\end{lemma}
\begin{proof}
	It suffices to prove \eqref{stab_fine} since one may apply same strategy to prove \eqref{stab_ms}.
	Take $q=\text{div}(\bu_{f})$ in \eqref{eq:fine sol p} and one have
	\begin{align*}
		\|\text{div}(\bu_{f})\|_{L^2}^2\leq \int_{\Omega }f\text{div}(\bu_{f})\leq \|f\|_{L^2}\|\text{div}(\bu_{f})\|_{L^2},
	\end{align*}
	where we apply the Cauchy-Schwartz inequality. Using the Poincar$\acute{e}$ inequality and $\bu_{f}\in V_f$, one can obtain
	\begin{align*}
		\|\bu_{f}\|_{L^2}\leq C\|\text{div}(\bu_{f})\|_{L^2}.
	\end{align*}
	Hence \eqref{stab_fine} is proved. Similarly, we can prove \eqref{stab_ms}.
\end{proof}
We now give estimates $e_1$ and $e_5$.
\begin{theorem}
	Recall that $\bu_{f}^{i,2}:=\bu_{f}(\bx;\omega_i,f_2)$ and $\bu_{\text{ms}}^{i,2}:=\bu_{\text{ms}}(\bx;\omega_i,f_2)$ are fine-grid and multiscale velocity solutions to \eqref{eq:fine sol v}-\eqref{eq:fine sol p} and \eqref{eq:ms sol v}-\eqref{eq:ms sol p} with $\omega_i$ and $f_2$, for $i=1,2$. Moreover, $p_{f}^{i,2}:=p_{f}(\bx;\omega_i,f_2)$ and $p_{\text{ms}}^{i,2}:=p_{\text{ms}}(\bx;\omega_i,f_2)$ are corresponding fine-grid and multiscale pressure solutions with $\omega_i$ and $f_2$. Then the following estimations hold,
	\begin{align}
		\|\kappa_2^{-\frac{1}{2}}\bu_{f}^{2,2}-\kappa_1^{-\frac{1}{2}}\bu_{f}^{1,2}\|_{L^2}&\leq C \|\kappa_2^{-\frac{1}{2}}-\kappa_1^{-\frac{1}{2}}\|_{L^{\infty}}\|f_2\|_{L^2}, \label{e_1}\\
		\|\kappa_2^{-\frac{1}{2}}\bu_{\text{ms}}^{2,2}-\kappa_1^{-\frac{1}{2}}\bu_{\text{ms}}^{1,2}\|_{L^2}&\leq C\|\kappa_2^{-\frac{1}{2}}-\kappa_1^{-\frac{1}{2}}\|_{L^{\infty}}\|f_2\|_{L^2}. \label{e_5}
	\end{align}
	\label{thm:e_1}
\end{theorem} 
\begin{proof}
	Recalling \eqref{eq:fine sol v} and \eqref{eq:fine sol p} associated with $\kappa_2$ and $f_2$, we obtain
	\begin{align}
		\int_\Omega  \kappa_2^{-1}\bu_f^{2,2}\cdot\bv-\int_\Omega\text{div}(\bv)p_f^{2,2}&=0,\quad\forall \bv\in V_f^0,\label{1_1} \\ 
		\int_\Omega  \text{div}(\bu_f^{2,2})q&=\int_\Omega  f_2q,\quad\forall q\in Q_f.\label{1_2} 
	\end{align}
	Similarly, with $\kappa_1$ and $f_2$, one can attain
	\begin{align}
		\int_\Omega  \kappa_1^{-1}\bu_f^{1,2}\cdot\bv-\int_\Omega\text{div}(\bv)p_f^{1,2}&=0,\quad\forall \bv\in V_f^0, \\ 
		\int_\Omega  \text{div}(\bu_f^{1,2})q&=\int_\Omega  f_2q,\quad\forall q\in Q_f.\label{1_3}
	\end{align}
	Taking $\bv=\bu_f^{2,2}$ in \eqref{1_1} and $q=p_f^{2,2}$ in \eqref{1_2}, it holds that
	\begin{align}
		\int_\Omega  \kappa_2^{-1}|\bu_f^{2,2}|^2-\int_\Omega\text{div}(\bu_f^{2,2})p_f^{2,2}&=0, \\ 
		\int_\Omega  \text{div}(\bu_f^{2,2})p_f^{2,2}&=\int_\Omega  f_2p_f^{2,2}.
	\end{align}
	It follows that
	\begin{align}
		\int_\Omega  \kappa_2^{-1}|\bu_f^{2,2}|^2=\int_\Omega  f_2p_f^{2,2}. 
	\end{align}
	Taking $\bv=\bu_f^{1,2}$ in \eqref{1_1} and $q=p_f^{2,2}$ in \eqref{1_3}, it holds that
	\begin{align}
		\int_\Omega  \kappa_2^{-1}\bu_f^{2,2}\cdot\bu_f^{1,2}=\int_\Omega  f_2p_f^{2,2}. 
	\end{align}
	Then we have 
	\begin{align}
		\int_\Omega  \kappa_2^{-1}\bu_f^{2,2}\cdot\bu_f^{1,2}=\int_\Omega  \kappa_2^{-1}|\bu_f^{2,2}|^2.\label{1_4} 
	\end{align}
	Similarly, one can derive that 
	\begin{align}
		\int_\Omega  \kappa_1^{-1}\bu_f^{1,2}\cdot\bu_f^{2,2}=\int_\Omega  \kappa_1^{-1}|\bu_f^{1,2}|^2. \label{1_5} 
	\end{align}
	Based on \eqref{1_4} and \eqref{1_5}, we have
	\begin{align*}
		\int_{\Omega}|\kappa_1^{-\frac{1}{2}}\bu_f^{1,2}-\kappa_2^{-\frac{1}{2}}\bu_f^{2,2}|^2&=\int_{\Omega}|\kappa_1^{-\frac{1}{2}}\bu_f^{1,2}|^2+|\kappa_2^{-\frac{1}{2}}\bu_f^{2,2}|^2-2(\kappa_1\kappa_2)^{-\frac{1}{2}}\bu_f^{1,2}\cdot \bu_f^{2,2}\\
		&=\int_{\Omega}|\kappa_1^{-\frac{1}{2}}-\kappa_2^{-\frac{1}{2}}|^2\bu_f^{1,2}\cdot \bu_f^{2,2}.
	\end{align*}
	We conclude that 
	\begin{align*}
		\|\kappa_1^{-\frac{1}{2}}\bu_f^{1,2}-\kappa_2^{-\frac{1}{2}}\bu_f^{2,2}\|_{L^2}^2&\leq \|\kappa_1^{-\frac{1}{2}}-\kappa_2^{-\frac{1}{2}}\|_{L^{\infty}}^2\|\bu_f^{1,2}\|_{L^2}\|\bu_f^{2,2}\|_{L^2},\\
		&\leq C\|\kappa_1^{-\frac{1}{2}}-\kappa_2^{-\frac{1}{2}}\|_{L^{\infty}}^2\|f_2\|_{L^2}^2,
	\end{align*}
	where the last inequality holds based on Lemma \ref{stab}. Similarly, \eqref{e_5} can be obtained.
\end{proof}
To estimate $e_2$, we define a weighted $L^2$ norm by $\|\bv\|_{\kappa^{-1}, \Omega}^2:=\int_{\Omega} \kappa^{-1}|u|^2$.
\begin{theorem}
	Recall that for $i=1,2$, $\bu_{f}^{1,i}:=\bu_{f}(\bx;\omega_1,f_i)$ and $\bu_{\text{ms}}^{1,i}:=\bu_{\text{ms}}(\bx;\omega_1,f_i)$ are fine-grid and multiscale velocity solutions to \eqref{eq:fine sol v}-\eqref{eq:fine sol p} and \eqref{eq:ms sol v}-\eqref{eq:ms sol p} with $\omega_1$ and $f_i$; $p_{f}^{1,i}:=p_{f}(\bx;\omega_1,f_i)$ and $p_{\text{ms}}^{1,i}:=p_{\text{ms}}(\bx;\omega_1,f_i)$ are corresponding fine-grid and multiscale pressure solutions with $\omega_1$ and $f_i$. Then, it holds that,
	\begin{align}
		\|\bu_{f}^{1,1}-\bu_{f}^{1,2}\|_{\kappa_1^{-1},\Omega}&\lesssim \kappa_{1,\text{min}}^{-\frac{1}{2}}\|f_1-f_2\|_{L^2},\label{e2_2}\\
		\|\bu_{\text{ms}}^{1,1}-\bu_{\text{ms}}^{1,2}\|_{\kappa_1^{-1},\Omega}&\lesssim C_{\text{infsup }}\kappa_{1,\text{min}}^{-\frac{1}{2}}\|f_1-f_2\|_{L^2},\label{e4}
	\end{align}
	where $\kappa_{1,\text{min}}:=\displaystyle\min_{\bx\in \Omega}\kappa(\bx;\omega_1)$.
\end{theorem}
\begin{proof}
	Recall that $\kappa_1:=\kappa(\bx;\omega_1)$.  Similarly as the proof in Theorem \ref{thm:e_1}, we subtract \eqref{eq:fine sol v} and \eqref{eq:fine sol p} corresponding with $f_1$ and $f_2$ to get 
	\begin{align}
		\int_{\Omega} \kappa_1^{-1}(\bu_{f}^{1,1}-\bu_{f}^{1,2})\cdot \bv-\int_{\Omega} \text{div}(\bv)(p_{f}^{1,1}-p_{f}^{1,2})&=0, \quad\forall \bv\in V_f^0, \label{middle_3}\\ 
		\int_{\Omega} \text{div}(\bu_{f}^{1,1}-\bu_{f}^{1,2})q&=\int_{\Omega} (f_1-f_2)q, \quad\forall q\in Q_f.
	\end{align}
	Take $\bv=\bu_{f}^{1,1}-\bu_{f}^{1,2}$, $q=p_{f}^{1,1}-p_{f}^{1,2}$ and it holds that
	\begin{align*}
		\|\bu_{f}^{1,1}-\bu_{f}^{1,2}\|_{\kappa_1^{-1},\Omega}^2=\int_{\Omega}(f_1-f_2)(p_{f}^{1,1}-p_{f}^{1,2}).
	\end{align*}
	By Cauchy-Schwartz inequality, 
	\begin{align}
		\|\bu_{f}^{1,1}-\bu_{f}^{1,2}\|_{\kappa^{-1},\Omega}^2\leq \|p_{f}^{1,1}-p_{f}^{1,2}\|_{L^2}\|f_1-f_2\|_{L^2}. \label{middle_4}
	\end{align}
	The Raviart-Thomas elements satisfy the following inf-sup condition \cite{brezzi1991variational}:
	\begin{align}
		\|q_h\|_{L^2}\lesssim \sup_{\bv_h\in V_f} \frac{\int_{\Omega}\text{div}(\bv_h)q_h}{\|v_h\|_{H(\text{div};\Omega)}},
		\quad \forall q_h\in Q_h.
	\end{align}
	Using \eqref{middle_3}, we obtain
	\begin{align}
		\|p_{f}^{1,1}-p_{f}^{1,2}\|_{L^2}\lesssim \kappa_{1,\text{min}}^{-\frac{1}{2}}\|\bu_{f}^{1,1}-\bu_{f}^{1,2}\|_{\kappa_1^{-1},\Omega}.\label{middle_5}
	\end{align}
	By \eqref{middle_4} and \eqref{middle_5}, we conclude that
	\begin{align}
		\|\bu_{f}^{1,1}-\bu_{f}^{1,2}\|_{\kappa_1^{-1},\Omega}\lesssim \kappa_{1,\text{min}}^{-\frac{1}{2}}\|f_1-f_2\|_{L^2}.
	\end{align}
	The process to prove \eqref{e4} is similar. We subtract \eqref{eq:ms sol v} and \eqref{eq:ms sol p} associated with $f_1$ and $f_2$ and it holds that
	\begin{align}
		\int_{\Omega} \kappa_1^{-1}(\bu_{\text{ms}}^{1,1}-\bu_{\text{ms}}^{1,2}) \cdot\bv-\int_{\Omega} \text{div}(\bv)(p_{\text{ms}}^{1,1}-p_{\text{ms}}^{1,2})&=0, \quad\forall \bv\in V_{\text{ms}}, \label{middle_6}\\ 
		\int_{\Omega} \text{div}(\bu_{\text{ms}}^{1,1}-\bu_{\text{ms}}^{1,2})q&=\int_{\Omega} (f_1-f_2)q, \quad\forall q\in Q_{\text{ms}}.
	\end{align}
	Taking $\bv=\bu_{\text{ms}}^{1,1}-\bu_{\text{ms}}^{1,2}$ and $q=p_{\text{ms}}^{1,1}-p_{\text{ms}}^{1,2}$, it follows that
	\begin{align}
		\|\bu_{\text{ms}}^{1,1}-\bu_{\text{ms}}^{1,2}\|_{\kappa_1^{-1},\Omega}^2\leq \|p_{\text{ms}}^{1,1}-p_{\text{ms}}^{1,2}\|_{L^2}\|f_1-f_2\|_{L^2}. \label{middle_7}
	\end{align}
	From Theorem 4.2 in \cite{chung2015mixed}, we have
	\begin{align}
		\|p\|_{L^2}\lesssim C_{\text{infsup}}\sup_{\bv\in V_{\text{ms},0}} \frac{\int_{\Omega}\text{div}(\bv)p}{\|v_h\|_{H(\text{div};\Omega_1^{-1})}},
		\quad \forall p\in Q_{\text{ms}},
	\end{align}
	where $V_{\text{ms},0}$ is the offline multiscale space.
	\begin{align}
		C_{\text{infsup}}=(\displaystyle 1+\max_{1\leq i\leq N_{\text{E,c}}}\min_r\int_{D_i}\kappa_1^{-1}|\psi_r^{i}|^2)^{\frac{1}{2}}, \label{c_infsup}
	\end{align}  
	where the minimum is taken over all indices $r$ with the property $\int_{E_i} \psi_r^{i}\cdot m_i\neq 0$. Here $m_i$ is an outer normal vector to $E_i$.
	Since $V_{\text{ms}}$ is an enriched multiscale space from $V_{\text{ms},0}$, $V_{\text{ms},0}\subset V_{\text{ms}}$. Then the following holds,
	\begin{align}
		\|p\|_{L^2}\lesssim C_{\text{infsup}}\sup_{\bv\in V_{\text{ms}}} \frac{\int_{\Omega}\text{div}(\bv)p}{\|v_h\|_{H(\text{div};\Omega_1^{-1})}},
		\quad \forall p\in Q_{\text{ms}},
	\end{align}
	Using \eqref{middle_6}, we have
	\begin{align}
		\|p_{\text{ms}}^{1,1}-p_{\text{ms}}^{1,2}\|_{L^2}\lesssim \kappa_{1,\text{min}}^{-\frac{1}{2}}\|\bu_{\text{ms}}^{1,1}-\bu_{\text{ms}}^{1,2}\|_{\kappa_1^{-1},\Omega}.\label{middle_8}
	\end{align}
	
	Combining \eqref{middle_7}-\eqref{middle_8}, \eqref{e4} is attained.
\end{proof}
In the following part, we estimate $e_3$. 
Given a region $D$, we recall that $V_D:=\bigoplus_{D_i\subset D}V_{\text{snap}}^{i}$. $\tilde{V}_D$ is defined to be the divergence free subspace of $V_D$.
Then we define the residual operator $R_D^{(m)}$ as a linear functional on $V_{D}$ on the enrichment level $m$ by
\begin{align*}
	R_{D}^{(m)}(v)=\int_{D}\kappa_1^{-1} \bu_{\text{ms}}^{(m)}\cdot v-\int_{D}\text{div}(v)p_{\text{ms}}^{(m)}, \quad v\in V_{D}.
\end{align*}
If we restrict $R_D^{(m)}$ on $\tilV_D$, it holds that 
\begin{align*}
	R_{D}^{(m)}(v)=\int_{D}\kappa_1^{-1} \bu_{\text{ms}}^{(m)}\cdot v.
\end{align*}
Moreover, we define the norm of $R_D^{(m)}$ by
\begin{align*}
	\|R_D^{(m)}\|_{({V_{D}})^{*}}=\sup_{v\in V_{D}}\dfrac{|R_D^{(m)}(v)|}{\|v\|_{H(\text{div};D;\kappa_1^{-1})}}.
\end{align*}

\begin{theorem}\cite{online-mixed}
	Recall that $\bu_{f}^{1,1}:=\bu_{f}(\bx;\omega_1,f_1)$, $\bu_{\text{ms}}^{1,1}:=\bu_{\text{ms}}(\bx;\omega_1,f_1)$ are fine-grid and multiscale velocity solutions to \eqref{eq:fine sol v}-\eqref{eq:fine sol p} and \eqref{eq:ms sol v}-\eqref{eq:ms sol p} with $\omega_1$ and $f_1$. Suppose $\bu_{\text{ms}}^{1,1}$ is the multiscale solution at the enrichment level $m$, i.e., $\bu_{\text{ms}}^{1,1}=\bu_{\text{ms}}^{1,1,(m)}$, which is the solution to \eqref{ms_k}  with $\omega_1$ and $f_1$. Then we have the following estimation:
	\begin{align*}
		\|\bu_{f}^{1,1}-\bu_{\text{ms}}^{1,1}\|_{\kappa_1^{-1},\Omega}\leq C_{\text{err}}\sum_{j=1}^{Ne}\|R_{D_j}^{(0)}\|_{(V_{\text{snap}}^{j})^{*}}(\lambda_{l_j+1}^{(j)})^{-1}-
		\sum_{i=1}^{m}\sum_{j=1}^{P}\|R_{D_j}^{(i)}\|_{(\hat{V}_{D_j})^{*}},
	\end{align*}
	where $C_{\text{err}}=\frac{C_V H}{h}$ and $C_V$ depends on the polynomial order of the fine-grid basis functions in $V_{\text{snap}}$.
\end{theorem}

\begin{theorem}
	For each $\omega_2$ and $f_2$, we recall that $\bu_{f}^{2,2}:=\bu_{f}(\bx;\omega_2,f_2)$, $\bu_{\text{ms}}^{2,2}:=\bu_{\text{ms}}(\bx;\omega_2,f_2)$ are fine-grid and multiscale velocity solutions to \eqref{eq:fine sol v}-\eqref{eq:fine sol p} and \eqref{eq:ms sol v}-\eqref{eq:ms sol p} with $\omega_2$ and $f_2$. Combining the previous conclusions, the following estimation holds,
	\begin{eqnarray}
		\begin{aligned}
		\|\bu_{f}^{2,2}-\bu_{\text{ms}}^{2,2}\|_{\kappa_2^{-1},\Omega}&\lesssim C \|\kappa_2^{-\frac{1}{2}}-\kappa_1^{-\frac{1}{2}}\|_{L^{\infty}}\|f_2\|_{L^2}+
(1+C_{\text{infsup }})\kappa_{1,\text{min}}^{-\frac{1}{2}}\|f_1-f_2\|_{L^2}+\\
&C_{\text{err}}\sum_{j=1}^{Ne}\|R_{D_j}^{(0)}\|_{(V_{\text{snap}}^{j})^{*}}(\lambda_{l_j+1}^{(j)})^{-1}-
\sum_{i=1}^{m}\sum_{j=1}^{P}\|R_{D_j}^{(i)}\|_{(\hat{V}_{D_j})^{*}}.
		\end{aligned}
	\end{eqnarray}
Here we recall that $C_{\text{infsup }}$ is defined  in \eqref{c_infsup}.
\end{theorem}

We remark that the this theorem states that the total error is a consequence of three parts, i.e. different sources $f_1$ and $f_2$, different permeability fields $\kappa_1$ and $\kappa_2$, and multiscale approximation error, respectively. As for the last part, it can be interpreted as a truncated error since only a few eigenvectors of the local spectral problems are included in the offline  stage I. However, since there is a second stage, the offline stage II, we can compensate the previous approximation error by some residual-driven bases, which is characterized by subtracting the summation of some residuals. 

\section{Numerical results}
In this part, we present numerical results. 
In the following tests, we assume that the logarithmic permeability field $Y(\bx,\omega)=\log(\kappa(\bx,\omega))$ is a second-order stationary Gaussian random field. To perform the KL expansion of the random field, we need to define a covariance function as follows:
\begin{align}
	C(\bx,\bz)=\sigma^2\exp\left(-\frac{|x_1-z_1|^2}{2\eta_1^2}-\frac{|x_2-z_2|^2}{2\eta_2^2}\right), \label{cov_fun}
\end{align}
where $\bx=(x_1,x_2)$ and $\bz=(z_1,z_2)$ are two arbitrary variables in the spatial domain $\Omega$. Moreover, $\eta_1$ and $\eta_2$ are correlation lengths in first and second directions. Suppose $(\lambda_i,f_i)$, $i=1,\ldots,+\infty$, are eigenvalues and corresponding eigenfunctions of $C(\bx,\bz)$. Without loss of generality, $\{\lambda_i\}$ are arranged in a nonincreasing order. One can approximate $Y(\bx,\omega)$ by 
\begin{align}
	Y(\bx,\omega)\approx E[Y(\bx,\omega)]+\sum_{i=1}^{N_k}\mu_i\sqrt{\lambda_i}f_i.
\end{align} 
where $\{\mu_i\}_{i=1}^{N_k}$ are standard identically independent Gaussian random variables. 

We consider five deterministic fields to serve as $E[Y(\bx,\omega)]$, denoted by $\kappa_i$, $i=1,\ldots,5$. Moreover, we set $\eta_1=\eta_2=\eta$.
We both apply two-dimensional and three-dimensional models. In Table \ref{meshsize}, we show the number of truncated KL expansion terms $N_k$ used in five models with two different choices of correlation lengths $\eta$. More specifically, we use $\eta=\frac{1}{8},\frac{1}{16}$ in three 2D permeability fields and for 3D case, $\eta=\frac{1}{4},\frac{1}{8}$. Here, we emphasize that due to a fast decay of eigenvalues of $C(\bx,\bz)$, the summation of a small number of largest eigenvalues can account for almost all the energy. Thus, the eigenvectors corresponding to the dominating eigenvalues can be used to represent the Gaussian field without losing much important information. Thanks to to this design, the computation effort contained in parameterizing the stochastic permeability field is greatly reduced, which is another contribution of our proposed reduced-order method.

$\Omega$ is set to be $[0,2.2]\times[0,0.6]$ for $\kappa_1$ and $\kappa_2$. For $\kappa_3$, $\Omega=[0,1.28]\times[0,1.28]$. As for the three-dimensional case, $\Omega=[0,2.2]\times [0,0.6] \times [0,0.3]$ for $\kappa_4$ and $[0,0.64]\times[0,0.64]$ for $\kappa_5$. $\kappa_1$, $\kappa_2$ and $\kappa_4$ are extracted from SPE10 model, which serves as a benchmark in estimating the effect of upscaling and multiscale (ms for short) approaches. The SPE model is characterized by a high contrast and strong heterogeneity. $\kappa_1$ is the 36-th layer and $\kappa_2$ are the 85-th layer. $\kappa_4$ is composed of the last 30 layers of SPE model. The $\kappa_3$ and $\kappa_5$ are two high-contrast models, where the contrast is as big as $10^4$. In Table \ref{meshsize}, we calculate the mesh sizes of fine-scale and coarse-scale mesh applied to five permeability fields. Before the discussion, we define $\bv_f$ and $\bv_{\text{ms}}$ to be the fine-grid solution  and the GMsFEM solution respectively. Furthermore, the relative velocity error $e_{\bv}$ is defined as follows,
\begin{align*}
	e_{\bv}:=\frac{\int_{\Omega}\kappa^{-1}|\bv_f-\bv_{\text{ms}}|^2}{\int_{\Omega}\kappa^{-1}|\bv_f|^2}.
\end{align*}
\begin{figure}[htbp!]
	\centering
	\subfigure[$\kappa_1$ in $\log_{10}$ scale ]{\includegraphics[trim={0cm 2cm 0cm 0cm}
		,clip,width=2.5in]{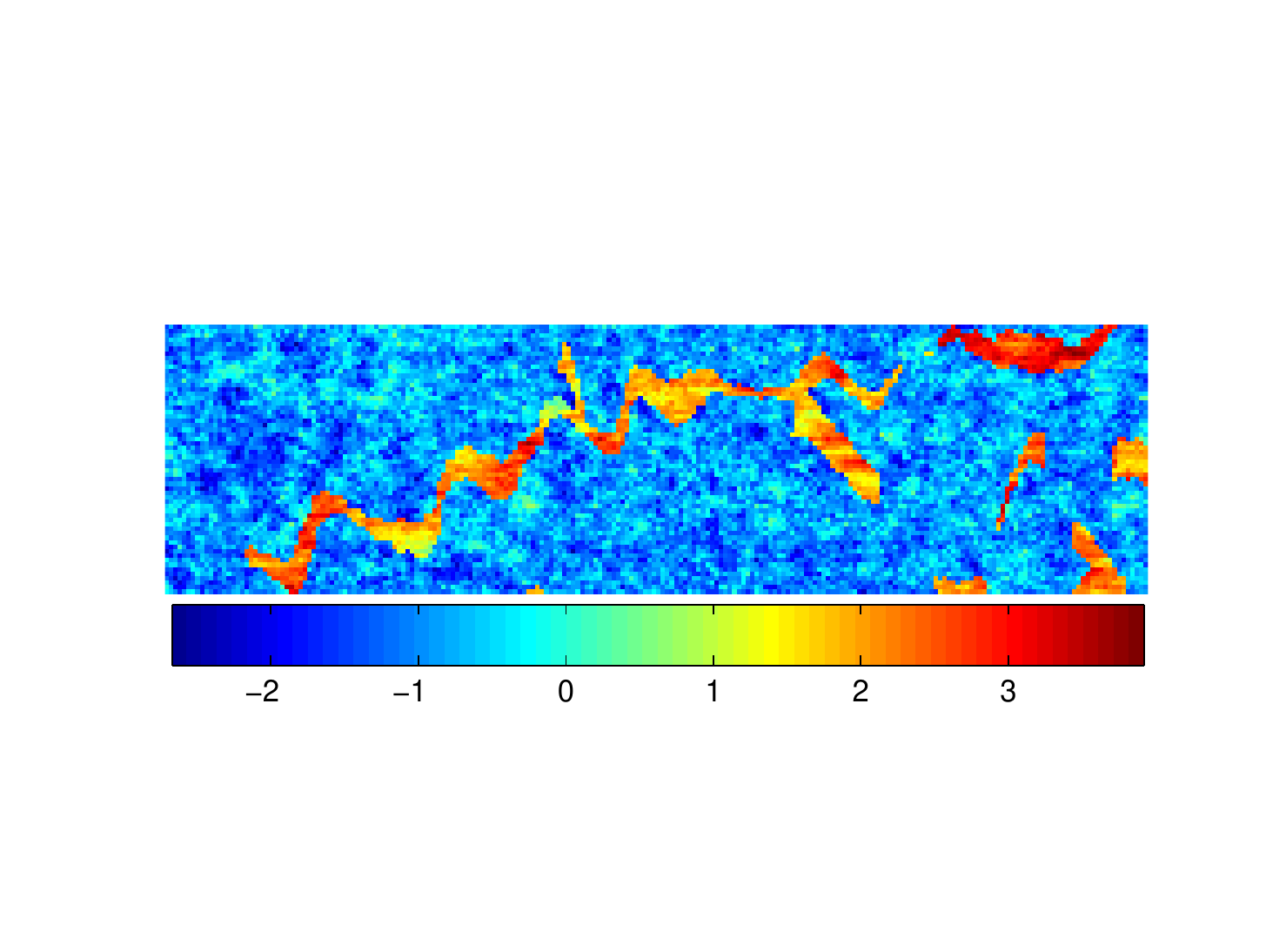}}
	\subfigure[$\kappa_2$ in $\log_{10}$ scale ]{\includegraphics[trim={0cm 2cm 0cm 0cm}
		,clip,width=2.5in]{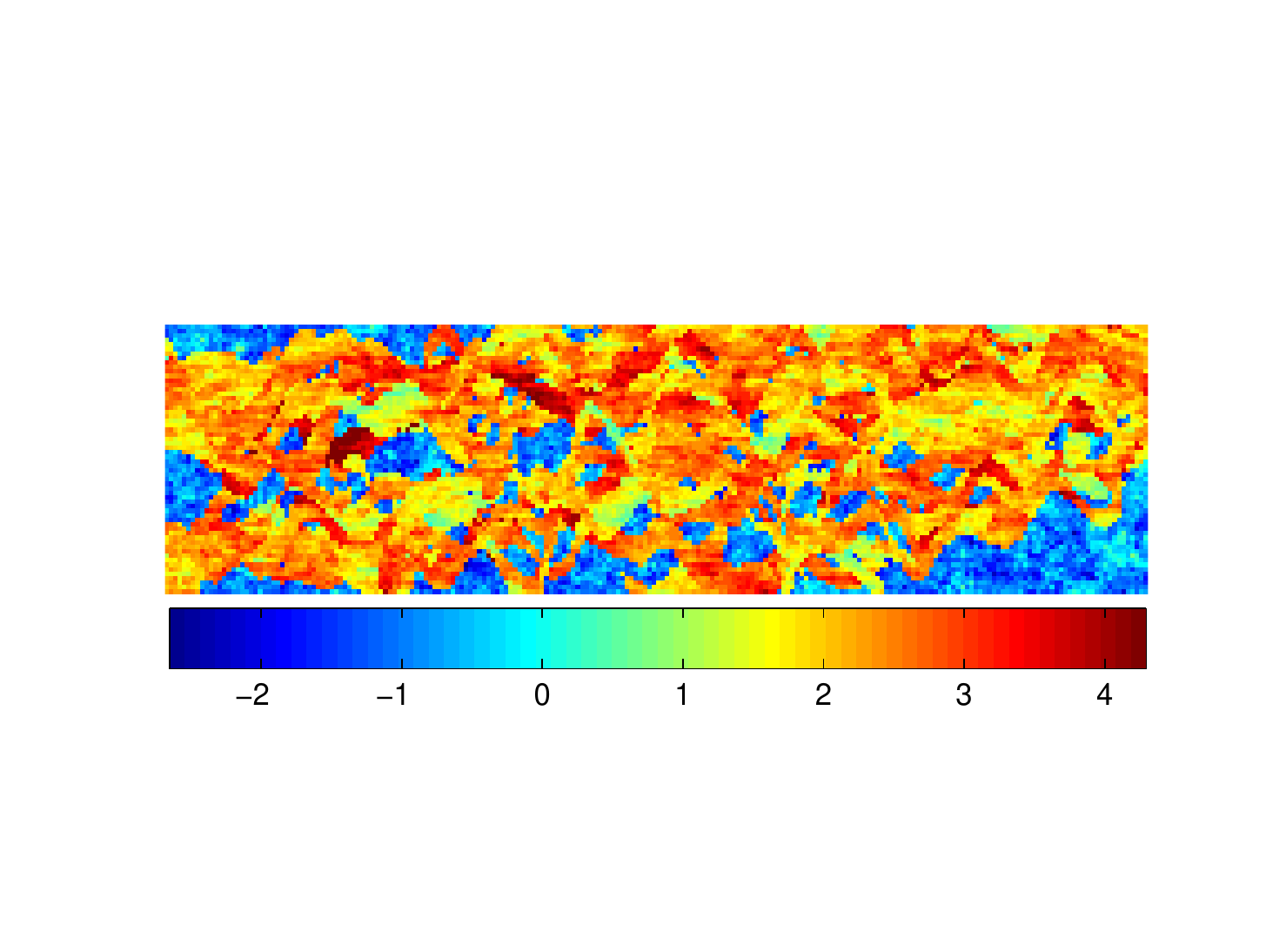}}	
	\subfigure[$\kappa_3$  ]{\includegraphics[trim={0cm 0cm 0cm 0cm}
		,clip,width=2.5in]{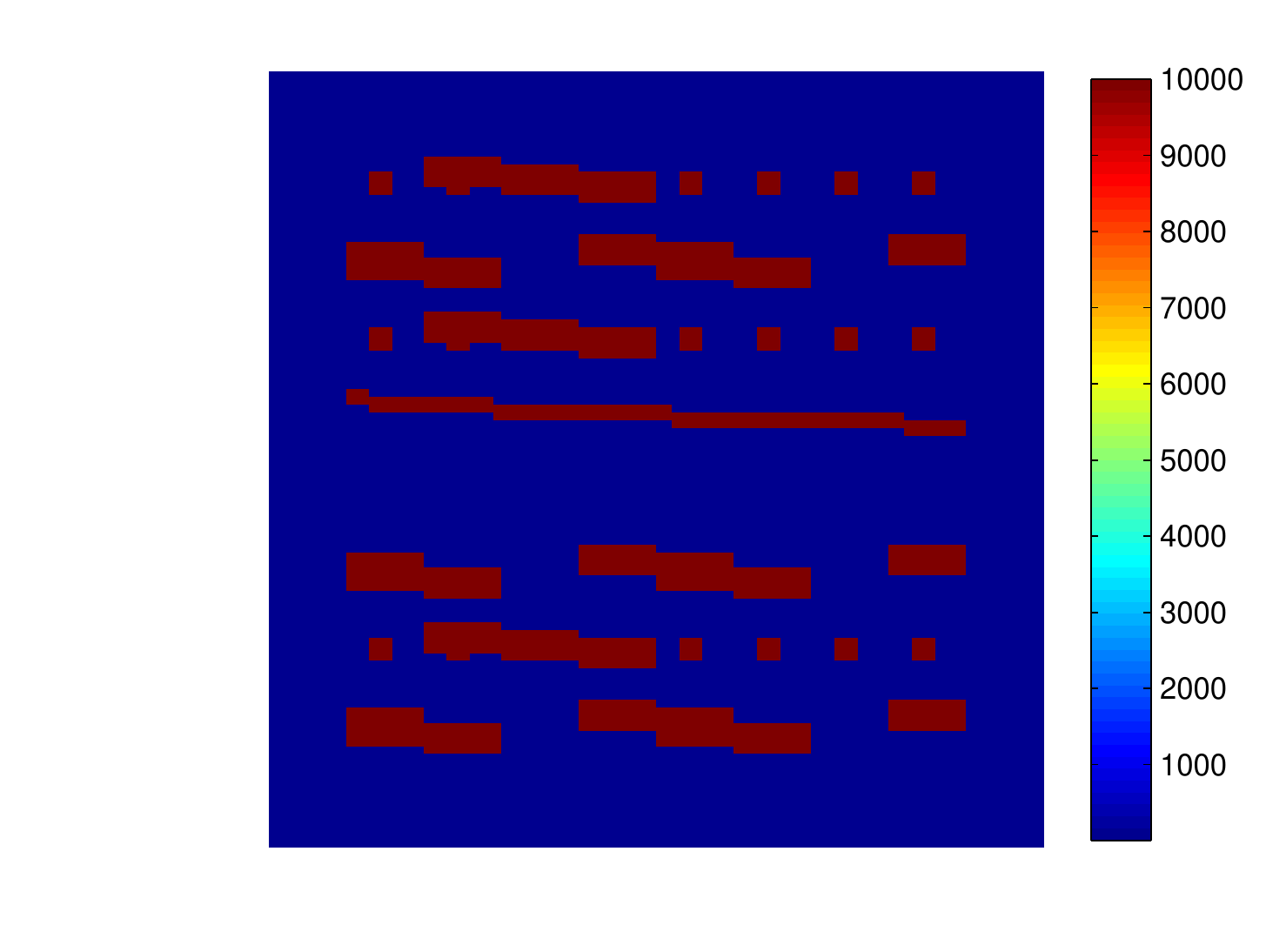}}
	\subfigure[$\kappa_4$ in $\log_{10}$ scale ]{\includegraphics[trim={0cm 0cm 0cm 0cm}
		,clip,width=2.5in]{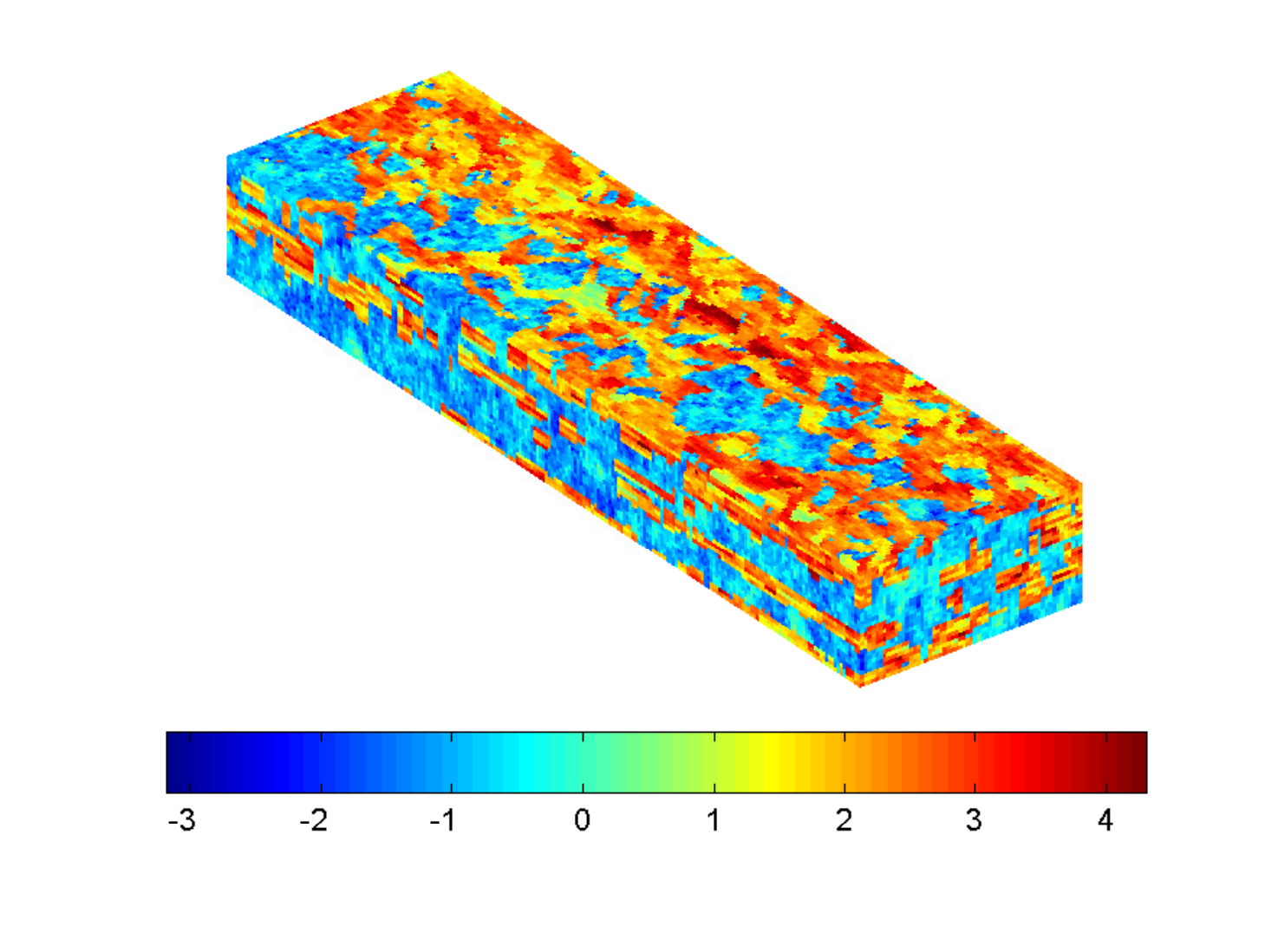}\label{spe10last30}}	
	\subfigure[$\kappa_5$]{\includegraphics[trim={3cm 1cm 3cm .5cm},clip,width=2.5in]{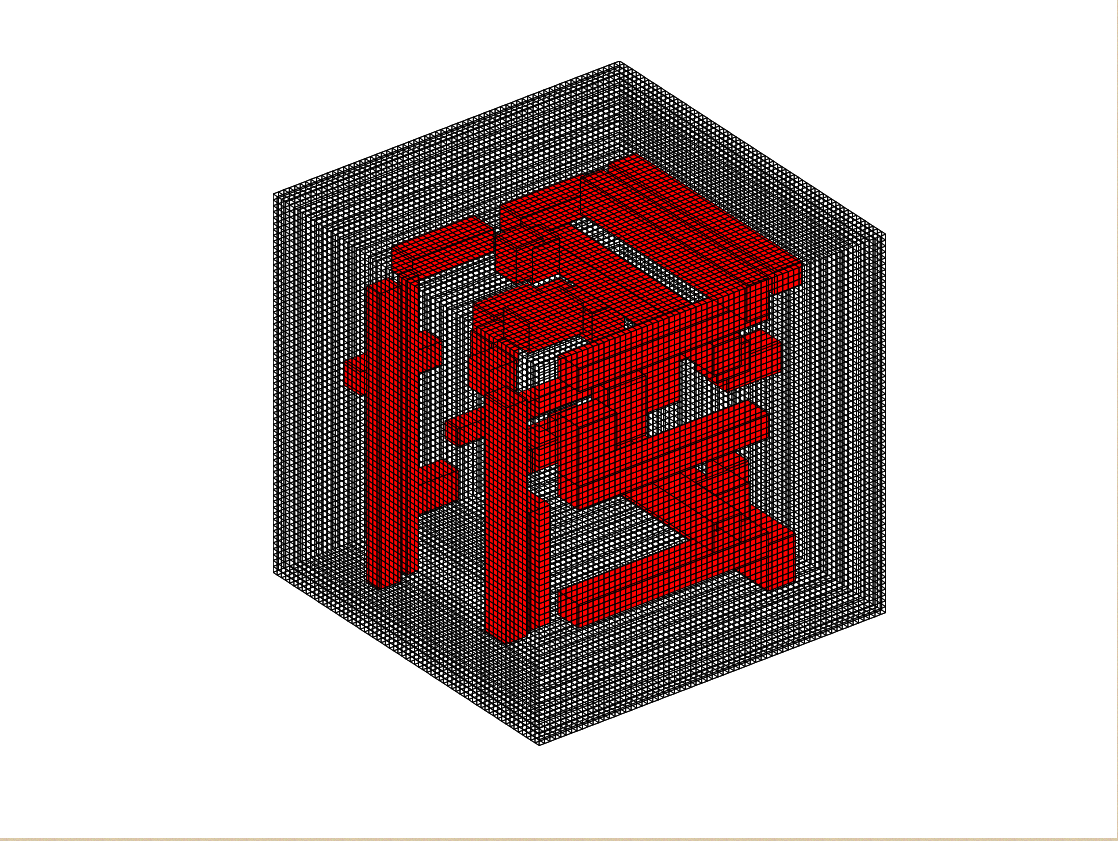}}	
	
	\caption{Permeability fields}
	\label{fig:models}
\end{figure}
We define $f_1$ and $f_2$ for 2D sources as follows.
\begin{align}
	f_1=\begin{cases}
		1, &x\in \zeta_1,\\
		-1, & x\in \zeta_4,\\
		0,  & otherwise.
	\end{cases}
	\quad
	f_2=\begin{cases}
		1, &x\in \zeta_i, i=1,\ldots,4\\
		-4, & x\in \zeta_5,\\
		0,  & otherwise,
	\end{cases}
	\label{source}
\end{align}
where $\zeta_i$, $i=1,\ldots,4$ are fine-scale elements at the four corners of $\Omega$. In particular, $\zeta_1=[0,h]\times [0,h]$, $\zeta_2=[L_x-h,L_x]\times [0,h]$, $\zeta_3=[0,h]\times [L_y-h,L_y]$, $\zeta_4=[L_x-h,L_x]\times[L_y-h,L_y]$. Moreover $\zeta_5$
is the fine-scale element at the center of $\Omega$. More specifically, $\zeta_5=[\frac{L_x-h}{2},\frac{L_x+h}{2}]\times [\frac{L_y-h}{2},\frac{L_y+h}{2}]$. We call $f_1$ a two-point source while $f_2$ is called five-point source. In 3D case, we can similarly define a two-point source and a five-point source as follows.
\begin{align}
	f_3=\begin{cases}
		1, &x\in \iota_1,\\
		-1, & x\in \iota_4,\\
		0,  & otherwise.
	\end{cases}
	\quad
	f_4=\begin{cases}
		1, &x\in \iota_i, i=1,\ldots,4\\
		-4, & x\in \iota_5,\\
		0,  & otherwise,
	\end{cases}
	\label{source_3d}
\end{align} 
where $\iota_i$ is the cube with bottom side $\zeta_i$, for $i=1,\ldots, 5$. In particular, $\iota_i=\zeta_i\times[0,L_z]$.  Unless specified, sources are chosen to be $f_1$ in 2D case and $f_3$ in 3D case. 

Before the presentation of numerical results, we first show a comparison of numbers of unknowns in solving a fine-scale flow equation and using multiscale bases to solve a flow equation. Without loss of generality, we can use a $2D$ case for illustration. Define $N_{\text{fine,x}}:=L_x/h$ and $N_{\text{fine,y}}:=L_y/h$ to be the numbers of fine-scale elements in the horizontal and vertical directions. Correspondingly, define $N_{\text{ms,x}}:=L_x/H$ and $N_{\text{ms,y}}:=L_y/H$ to be the numbers of coarse-scale  elements in two directions. Since th same number of muliscale bases are used in each local neighborhood, we let $N_{b}$ to be the dimension of each local multiscale space. Besides, we define $``A+B"$ to denote that A and B bases are incorporated in each local neighborhood in offline stage I and II, respectively. We then compare the DOF in fine-scale and coarse-scale equations. The Raviart-Thomas mixed finite element discretization of the flow equation gives rise to a linear system with $N_{\text{fine,x}}\times N_{\text{fine,y}}+(N_{\text{fine,x}}-1)\times N_{\text{fine,y}}+N_{\text{fine,x}}\times (N_{\text{fine,y}}-1)$ degrees of freedom.  In contrast, with multiscale bases, one can reduce the original flow equation to a system with $N_{\text{ms,x}}\times N_{\text{ms,y}}+N_{b}\times[(N_{\text{ms,x}}-1)\times N_{\text{ms,y}}+N_{\text{ms,x}}\times (N_{\text{ms,y}}-1)]$ degrees of freedom. In one example of our simulation, we choose $N_{\text{fine,x}}=220$ and $N_{\text{fine,y}}=60$, while $N_{\text{ms,x}}=11$ and $N_{\text{ms,y}}=3$. Besides, we choose $N_b=3$, which is sufficient to give a relatively good approximation in our cases. Under the above settings, the fine-scale equation system has $39320$ unknowns while the multiscale solver can result in a system with $189$ unknowns. From this example, we can see a great reduction is achieved by using multiscale bases. 

In 3D case we observe a more significant reduction. In Table \ref{k4 velo error} and \ref{k5 velo error}, we show Dof, training time as well as testing time for SPE model ($\kappa_4$) and high-contrast model ($\kappa_5$), respectively. In particular, the training time refers to the time spent in constructing multiscale bases while the testing time counts the time cost in solving a lower-dimensional system once based on the prepared multiscale bases. We first present the DOF and solving time of full-order system. For SPE model, the number of unknowns is 1562400 and a single solve needs 58.3s. As for the high-contrast model, the fine-grid system has 1036288 unknowns and it requires 31.86s per single solve. So it is time-consuming to solve such big systems. On the contrary, the multiscale systems are relatively smaller. Since $``2+1"$ bases contribute to a relatively accurate solution, we use this case for a comparison. For $\kappa_4$, $``2+1"$  bases give rise to a system with 3312 unknowns and the corresponding solving time is 2.41s per solve, which is remarkably shortened. For $\kappa_5$, the system with $2+1$ bases has DOF 4544 and each solve only needs about 2.92s. Even though the training process will increase computation time, about 90s with $\kappa_4$, the training process is only needed once. When there are numerous samples, this preparation effort is obviously cost-effective.

\subsection{Single-phase flow}
In this subsection, we consider a single-phase flow modeled by equations \eqref{model_v} and \eqref{model_p}.
In Table \ref{1/8 velo error2d}, we show the velocity errors corresponding to $\kappa_1$, $\kappa_2$ and $\kappa_3$ with different numbers of local bases when the correlation length $\eta=\frac{1}{8}$. More specifically, we both consider offline cases (from $``3+0"$ to $``20+0"$) and residual-driven cases (from $``1+1"$ to $``4+2"$). Here, we use the $``20+0"$ error to serve as a lower bound of multiscale approximation biases with $\kappa_1$ and $\kappa_2$ since in this case the velocity multiscale space is as big as the velocity snapshot space. We review that the dimension of a local velocity snapshot space is equal to the number of fine-scale edges on a single coarse-scale edge. From observations, we can obtain two main conclusions. Firstly, as more offline bases are used, the velocity error steadily decreases. As one can observe from the Table \ref{1/8 velo error2d}, the error of $\kappa_1$ decreases from $22.1\%$ at $``3+0"$ to $3.11\%$ at $``20+0"$, where 20 is the maximum for offline bases used in each local region. Secondly, the residual-driven bases can remarkably reduce errors. In other words, one can obtain a good accuracy with relatively few residual-driven bases. The error ($3.47\%$) corresponding to $``2+1"$ case in $\kappa_1$ is comparable to the $``20+0"$ case ($3.11\%$), which is sufficient to demonstrate the power of residual-driven bases. Since the $``2+1"$ case is nearly as accurate as the ``best'' case (``20+0''), it is reasonable that adding more residual-driven bases can not provide an apparent error reduction. As one can see, the ``2+1'' case has error $3.47\%$ while the $``2+2"$ has error $3.34\%$.  The error reduction from adding residual-driven bases is also significant when $\kappa_3$ is used. The error in $``3+1"$ case is $0.37\%$, less than $1/15$ of that in the  $``3+0"$ case. Apart from the two above conclusions, we can also observe that there is an evident accuracy improvement as one enriches the multiscale space from $``1+1"$ to $``2+1"$ bases. More specifically, errors associated with $\kappa_3$ decrease from $9.05\%$ to $0.88\%$ by starting with more than one offline bases.  The results from  $\eta=\frac{1}{16}$ are recorded in the Table \ref{1/16 velo error2d}, which share some similarities with the  case $\eta=\frac{1}{8}$. The $L_2$ errors corresponding to $\kappa_1$  decreases from $20.98\%$ at the $``3+0"$ case to $4.14\%$ at  the $``20+0"$ case. After adding one single residual-driven basis, the error straightly decreases to  $4.32\%$ at  the $``3+1"$ case, approximately as low as the $``20+0"$ case. Hence, only one residual-driven basis can improve the approximation such that it is nearly as accurate as the snapshot solution, which is sufficient to  show the strong power of residual-driven bases. Moreover, when we compare the above two cases ($\eta=\frac{1}{8}$, $\frac{1}{16}$), we can find that errors are slightly larger with a smaller correlation length $\eta$. This is resulted from the fact that when $\eta$ is larger, the covariance function defined in \eqref{cov_fun} will contribute to a more complicated permeability field while the velocity solutions are closely related with the permeability field. Hence, a smaller correlation length will increase the difficulty in approximating the velocity solutions.

\begin{figure}[htbp!]
	\centering
	\subfigure[$\eta=\frac{1}{8}$, different sources]{\includegraphics[width=0.32\textwidth]{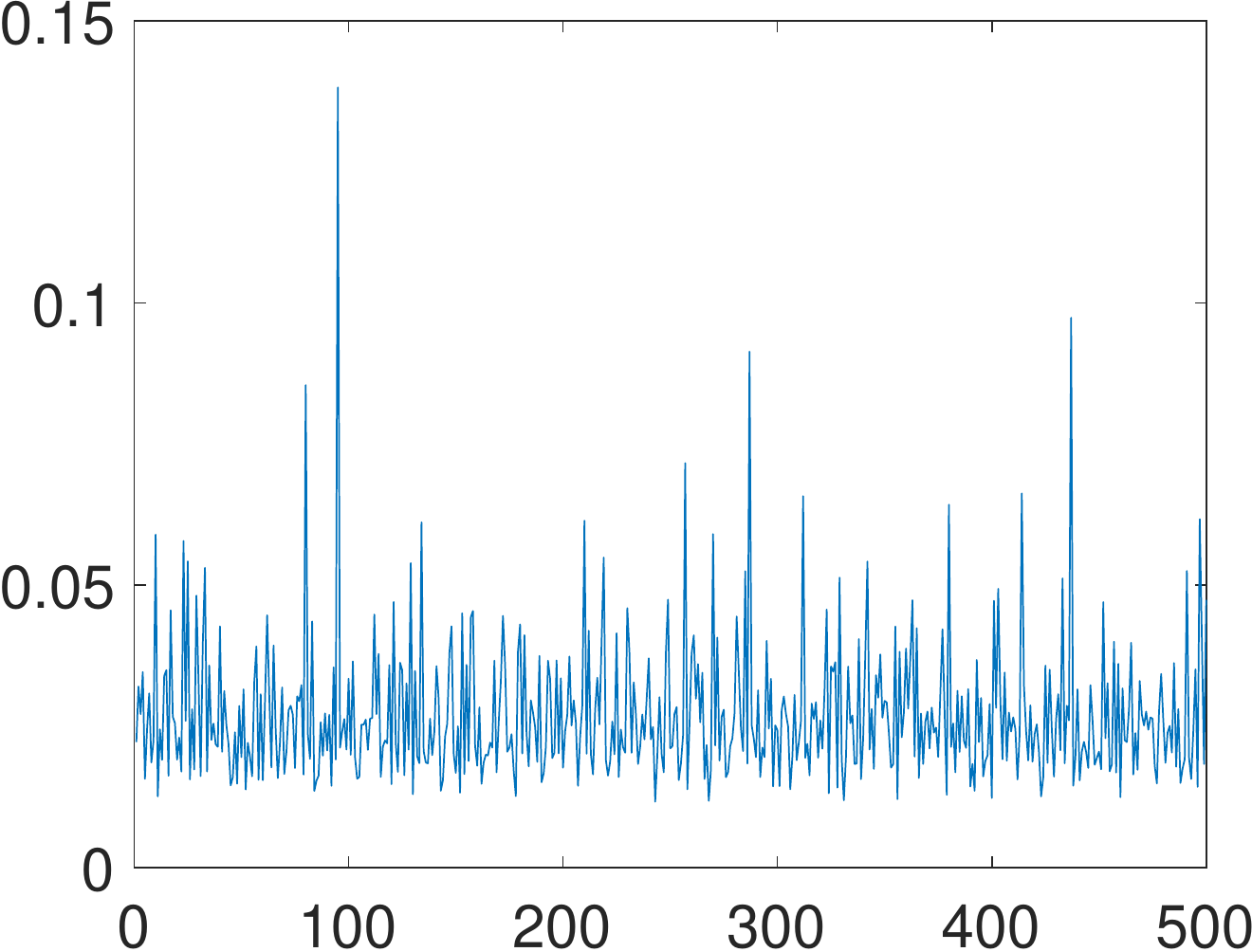}}
	\subfigure[$\eta=\frac{1}{8}$, same source]{\includegraphics[width=0.32\textwidth]{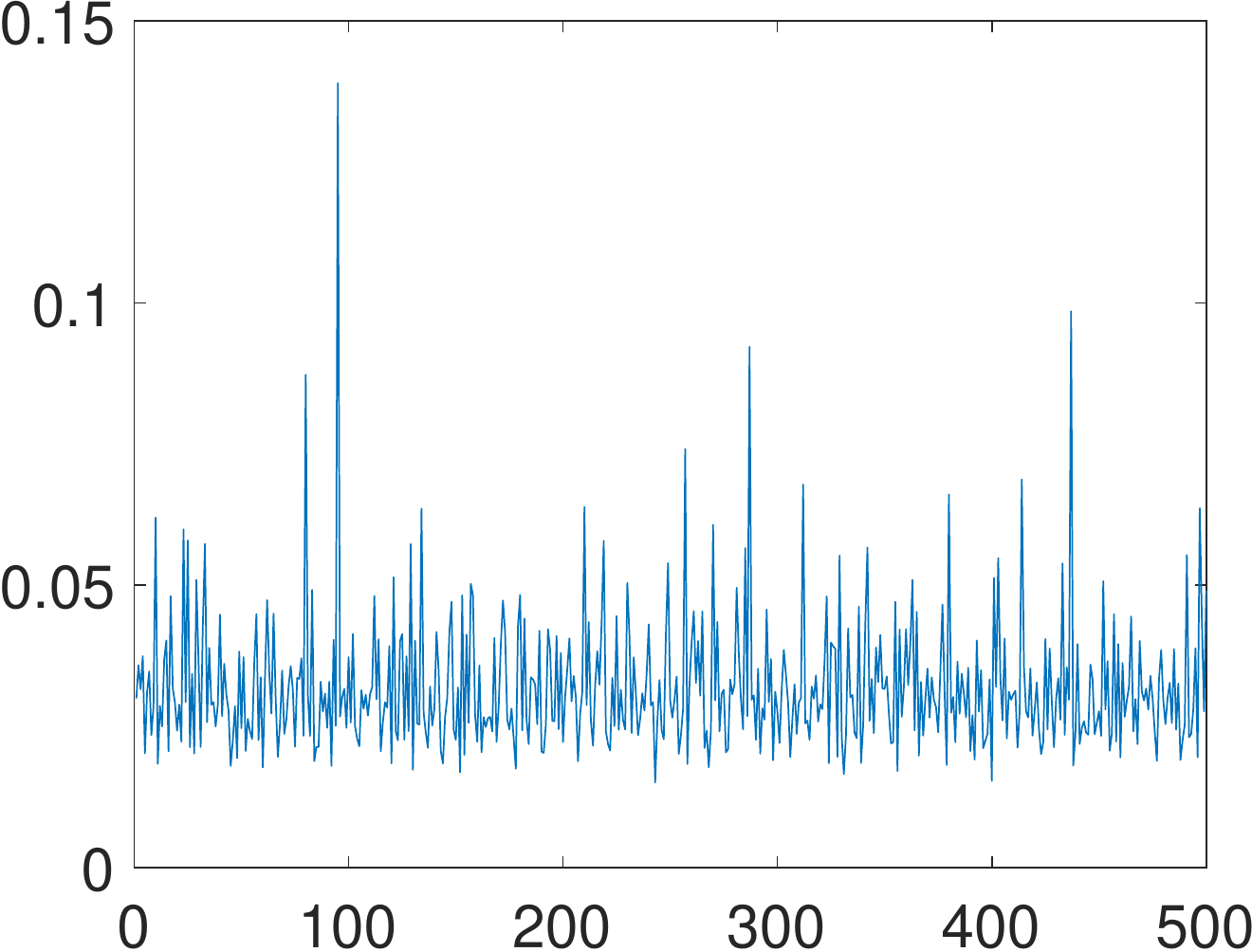}}
	\subfigure[$\eta=\frac{1}{16}$, same source]{\includegraphics[width=0.32\textwidth]{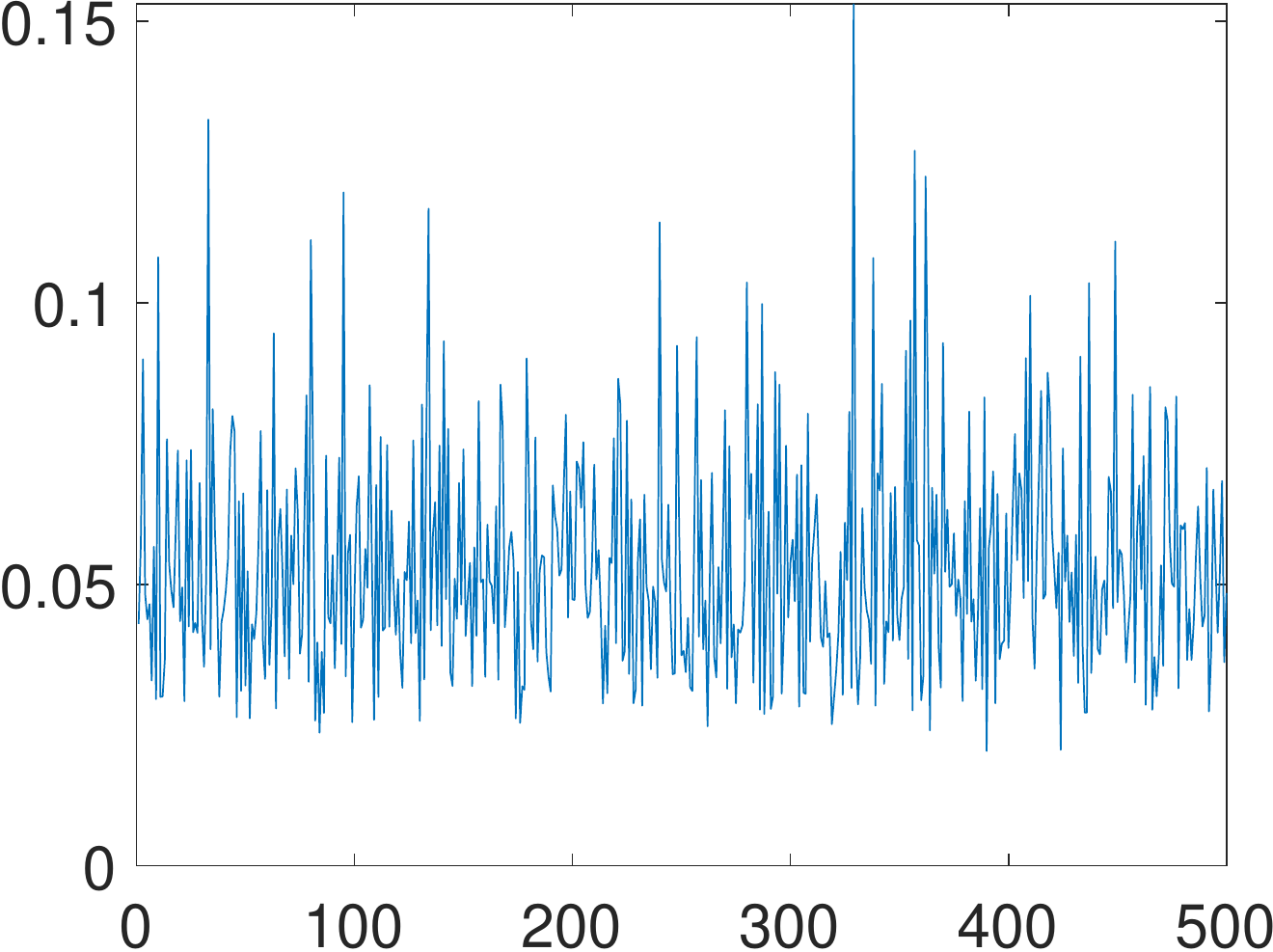}}
	\caption{Velocity errors of 500 samples with $\eta=1/8,1/16$. Left: $\eta=1/8$. Mean: 0.0320. Variance: $1.44\times 10^{-4}$.  We use $f_1$ in training while $f_2$ is used in testing. Middle: $\eta=1/8$. Mean: 0.0271. Variance: $1.59\times 10^{-4}$. $f_1$ is used in both training and testing. Right: $1/16$. Mean: 0.053. Variance: $3.98\times 10^{-4}$. $f_1$ is used in both training and testing. The number of multiscale bases are $``2+1"$. }
	\label{fig:vel_er}
\end{figure}

We now show the effect of correlation lengths on coarse-scale approximations. We further demonstrate the result when different sources are used in training and testing stages. In Figure \ref{fig:vel_er}, we plot velocity errors corresponding to 500 samples, where the correlation lengths are $\eta=\frac{1}{8},\frac{1}{16}$. For the smaller $\eta$, we both present the errors when the sources are different or the same in offline and online stages. With the same source $f_1$,  the mean error of $\eta=\frac{1}{8}$ is 0.0271 and the variance is $1.59\times 10^{-4}$. For $\eta=\frac{1}{16}$, the mean error is 0.052 and the variance is $3.98\times 10^{-4}$.  With $f_1$ in training and $f_2$ in testing, the mean error of  $\eta=\frac{1}{8}$ is 0.0320 and variance is $1.59\times 10^{-4}$. As one can verify from the figures, the errors corresponding to a smaller correlation length is more divergent than those corresponding to a larger correlation length. This reflects that there is a higher variation among different sample permeability fields that are generated with a smaller correlation length. Moreover, the errors corresponding to different sources in training and testing are only slightly larger than the ones with same sources, which to some extent shows a good generalization power of our method.

In the three-dimensional case, we choose $\eta=\frac{1}{4},\frac{1}{8}$ for both SPE model and high-contrast model. In the Table \ref{k4 velo error} and \ref{k5 velo error}, we show the errors associate with  $\kappa_4$ and $\kappa_5$, respectively. According to the numbers of fine-scale  elements in a single coarse-scale element, the dimensions of local velocity snapshot spaces for $\kappa_4$ and $\kappa_5$ are $100$ and $64$. Here, we use at most 32 offline multiscale bases in each local region. Since the numbers of offline bases have not attained the maximums, i.e. $100$ for $\kappa_4$ and $64$ for $\kappa_5$, the shown smallest velocity error without using residual-driven bases is attained with  $``32+0"$ bases, which is still relatively high, especially for $\kappa_4$ case with the error about $20\%$. Moreover, there is a remarkable accuracy improvement from the $``32+0"$ to the $``1+1"$ case with $\kappa_4$. In particular, with $\eta=\frac{1}{4}$,  the error $ 2\%$ in the $``1+1"$ case is less than $1/10$ of that in the $``3+0"$ case, which suffices to show the significance of residual-driven bases. On the other hand, the effect of adding residual-driven bases with $\kappa_5$ is not as significant as the one with $\kappa_4$  but it is still manifest. The error with $``1+1"$ bases accounts for nearly $1/8$ of that with $``32+0"$ bases. Hence, we can see that the error reduction in three-dimensional case is remarkable by including residual-driven bases. Also, the improvement is greater when the permeability field shows greater heterogeneity. Besides, we show the errors corresponding to a smaller correlation length $\frac{1}{8}$ in the last columns of the Table \ref{k4 velo error} and \ref{k5 velo error}. The error decay shares some similarities with the case $\eta=\frac{1}{4}$. One could observe that the residual-driven bases are more powerful in $\kappa_4$ than $\kappa_5$. Besides, due to a smaller correlation length, a higher level of randomness is included in the stochastic permeability field and thus yields a larger average error. When $\eta=\frac{1}{8}$, the smallest error $1\%$ corresponding to $\kappa_4$ is attained with the case of $``4+2"$, which is bigger than $0.5\%$ with $\eta=\frac{1}{4}$.

\begin{table}[htbp!]
	\centering
	\begin{tabular}{|c|c|c|c|c|c|}
		\hline
		& fine-scale mesh & coarse-scale mesh & $\eta=1/4$ &$\eta=1/8$ & $\eta=1/16$
		\tabularnewline\hline
		$\kappa_1$	& $220\times60$             & $11\times3$    	&      / & 38    & 136         \tabularnewline\hline
		$\kappa_2$  & $220\times60$             & $11\times3$       &     /  & 38    & 136              \tabularnewline\hline
		$\kappa_3$  & $128\times128$            & $8\times8$        &     /  & 46    &  150              \tabularnewline\hline
		$\kappa_4$  & $220\times60\times30$     & $22\times6\times3$  &     48          & 283   &   /           \tabularnewline\hline
		$\kappa_5$  & $64^3$                    & $8^3$               &     48         & 281   & /                \tabularnewline\hline
	\end{tabular}
	\caption{Fine-scale and coarse-scale mesh sizes (second and third columns), numbers of truncated KL expansion terms (last three columns) in 2D and 3D permeability fields with two different correlation lengths $\eta$. $\kappa_1$, $\kappa_2$ and $\kappa_3$ are three 2D permeability fields; $\kappa_4$ and $\kappa_5$ are two 3D permeability fields. In 2D case, we use $\eta=1/8,1/16$, while in 3D case, $\eta=1/4,1/8$.}
	\label{meshsize}
\end{table}

\begin{table}[H]
	\centering
	\begin{adjustbox}{max width=\textwidth}
		
		\begin{tabular}{|c|c|c|c|c|c|c|}
			\hline
			$N_b$&$e_{\text{v}}(\kappa_1)$&$e_{\text{v}}(\kappa_2)$ &$e_{\text{v}}(\kappa_3)$    \tabularnewline
			\hline
			3+0	&0.2110&0.0819 &0.0631 \tabularnewline\hline
			6+0	&0.1107&0.0426 &0.0335 \tabularnewline
			\hline
			8+0&0.0575&0.0379&  0.0302\tabularnewline\hline
			16+0&0.0317&0.0307&0.0036  \tabularnewline\hline
			20+0&0.0311&0.0305& / \tabularnewline\hline
			1+1	&0.0537&0.0540& 0.0905 \tabularnewline\hline
			2+1	&0.0347&0.0360&0.0088  \tabularnewline
			\hline
			2+2	&0.0334&0.0334&0.0039 \tabularnewline\hline
			3+1&0.0331&0.0330&0.0037 \tabularnewline\hline
			4+1&0.0323&0.0360&0.0037 \tabularnewline\hline
			5+1&0.0324&0.0317& 0.0037 \tabularnewline\hline
			4+2&0.0320&0.0319&0.0036 \tabularnewline\hline
		\end{tabular}
	\end{adjustbox}
	\caption{Average velocity errors corresponding to $\kappa_1$, $\kappa_2$ and $\kappa_3$ with the correlation length $1/8$. $500$ samples are included. }
	\label{1/8 velo error2d}
\end{table}

\begin{table}[H]
	\centering
	\begin{adjustbox}{max width=\textwidth}
		
		\begin{tabular}{|c|c|c|c|c|c|c|}
			\hline
			$N_b$&$e_{\text{v}}(\kappa_1)$&$e_{\text{v}}(\kappa_2)$ &$e_{\text{v}}(\kappa_3)$    \tabularnewline
			\hline
			3+0	&0.2098&0.0855&0.0640 \tabularnewline\hline
			6+0	&0.1130&0.0489&0.0349\tabularnewline
			\hline
			8+0&0.0626&0.0451& 0.0316\tabularnewline\hline
			16+0&0.0417&0.0383&0.0083  \tabularnewline\hline
			20+0&0.0414&0.0380&/  \tabularnewline\hline	
			1+1	&0.0619&0.0568& 0.1111 \tabularnewline\hline	
			2+1	&0.0441&0.0429&0.0116  \tabularnewline
			\hline
			2+2	&0.0435&0.0411&0.0087\tabularnewline\hline
			3+1& 0.0432&0.0406&0.0085 \tabularnewline\hline
			4+1&0.0431&0.0397&0.0085 \tabularnewline\hline
			5+1&0.0436& 0.0397&0.0084 \tabularnewline\hline
			4+2&0.0430&0.0398&0.0084 \tabularnewline\hline
		\end{tabular}
	\end{adjustbox}
	\caption{Average velocity errors corresponding to $\kappa_1$, $\kappa_2$ and $\kappa_3$ with the correlation length $1/16$. $500$ samples are included. }
	\label{1/16 velo error2d}
\end{table}
\begin{table}[H]
	\centering
	\begin{adjustbox}{max width=\textwidth}
		
		\begin{tabular}{|c|c|c|c|c|c|c|}
			\hline
			$N_b$& Dof&$T_{\text{train}}(s)$&$T_{\text{test}}(s)$&$e_{\text{v}}(\kappa_4;1/4)$ &$e_{\text{v}}(\kappa_4;1/8)$ \tabularnewline
			\hline
			3+0	&3312&10.56 & 2.54&0.256&0.257 \tabularnewline\hline
			8+0&8172 &19.30&6.07 &0.218&0.218\tabularnewline\hline
			16+0&15948& 38.71&15.98&0.212&0.212\tabularnewline\hline
			32+0&31500&90.97&61.90  &0.202&  0.201\tabularnewline\hline
			1+1	&2340 &96.05&2.11&0.020& 0.021\tabularnewline\hline			
			2+1	&3312 &91.13&2.41&0.011&0.013 \tabularnewline\hline
			2+2	&4284 &150.71&3.24&0.006&0.010\tabularnewline\hline
			2+3	&5256&211.45&3.95& 0.006&0.011 \tabularnewline\hline
			3+2	&5256 &98.21&3.24&0.006&0.011\tabularnewline\hline
			3+1 &4284&166.95&4.12& 0.008&0.011  \tabularnewline\hline
			4+1&5256 &111.79&3.97&0.008&0.012\tabularnewline\hline
			5+1&6228&122.00&4.54 &0.006&0.010 \tabularnewline\hline
			4+2&6228&187.16&4.49&0.005&0.010 \tabularnewline\hline
		\end{tabular}
	\end{adjustbox}
	\caption{Numbers of unknowns (Dof) and computational time for $\kappa_4$. Average velocity errors corresponding to $\kappa_4$ with correlation lengths $1/4$ and $1/8$. $500$ samples are included. $T_{\text{train}}$ is the CPU time for constructing the multiscale space and $T_{\text{test}}$ is the CPU time for solving a corresponding coarse-scale equation for a single sample.}
	\label{k4 velo error}
\end{table}

\begin{table}[H]
	\centering
	\begin{adjustbox}{max width=\textwidth}
		
		\begin{tabular}{|c|c|c|c|c|c|c|}
			\hline
			$N_b$& Dof&$T_{\text{train}}(s)$&$T_{\text{test}}(s)$&$e_{\text{v}}(\kappa_5;1/4)$ &$e_{\text{v}}(\kappa_5;1/8)$ \tabularnewline
			\hline
			3+0	&4544  &8.05 &2.91 &0.356&0.355  \tabularnewline\hline
			8+0 &11246  &22.01 &11.17  &0.235&0.235\tabularnewline\hline
			16+0 &22016 &64.65& 37.85&0.091& 0.091\tabularnewline\hline
			32+0 &43520 &216.70&149.94&0.031& 0.031\tabularnewline\hline
			1+1	 &3200 &49.22&1.50 &0.046&0.047 \tabularnewline\hline
			2+1	 &4544 &64.99&2.92&0.045&0.047 \tabularnewline\hline
			2+2	 &5888 &118.09&4.34&0.006&0.015 \tabularnewline\hline
			2+3	  &7232 &179.82&5.80&0.004&0.014\tabularnewline\hline
			3+1 &5888 &77.49&4.27& 0.044&0.045 \tabularnewline\hline
			3+2 &7232 &139.41&7.65	&0.006&0.014 \tabularnewline\hline
			4+1 &7232 &92.81&5.791& 0.043&0.045 \tabularnewline\hline
			4+2  &8576 &166.07  &7.62& 0.005&0.014\tabularnewline\hline
			5+1&8576 &107.04  &7.25&0.044& 0.045  \tabularnewline\hline
			
		\end{tabular}
	\end{adjustbox}
	\caption{ Numbers of unknowns (Dof) and computational time for $\kappa_5$. Average velocity errors corresponding to $\kappa_5$ with correlation lengths $1/4$ and $1/8$. $500$ samples are included. $T_{\text{train}}$ is the CPU time for constructing the multiscale space and $T_{\text{test}}$ is the CPU time for solving a corresponding coarse-scale equation for a single sample.}
	\label{k5 velo error}
\end{table}

\subsection{Two-phase flow and transport}
In this subsection, we consider a two-phase flow problem in a bounded domain $\Omega\in \mathbb{R}^2$. The two phases are water and oil, denoted by $w$ and $o$, respectively. In particular, water is injected to a heterogenous reservoir, which displaces the trapped oil to a production well.  To simplify the model, we assume there is no capillary pressure and only consider a gravity-free environment. Moreover, the two phases are assumed to be immiscible and incompressible. The two-phase problem is modeled by an elliptic equation coupled with a transport equation. Before presenting the equations, we first list some notations.
\begin{enumerate}
	\item $S$: water saturation; $\bv$: total Darcy velocity; $p$: pressure.
	\item $\kappa_{r\alpha}$, $\mu_{\alpha}$: relative permeability and viscosity of phase $\alpha$, $\alpha=w,o$.
	\item $f$, $r$ are external forcing terms.
	\item $\xi(S)$: total mobility.
	\item $F(S)$: fractional flow function for water phase.
	\item $\partial D_{\text{out}}$: outer flow boundary.
\end{enumerate}
Under the above assumptions, the two-phase flow problem is formulated as a simplified coupling system, where the flow equation is as follows.
\begin{eqnarray}
	\begin{aligned}
		-\xi(S)\kappa(\bx;\omega)\nabla p&=\bv \text{ in }\Omega,\\
		\text{div}(\bv) &=f \text{ in }\Omega,\\
		\bv\cdot \bn&=0\text{ on }\partial \Omega.
		\label{flow_eqn}
	\end{aligned}
\end{eqnarray}
$\xi(S)$ is defined by
\begin{align*}
	\xi(S):=\frac{\kappa_{rw}(S)}{\mu_{w}}+\frac{\kappa_{ro}(S)}{\mu_{o}}.
\end{align*}
Here
\begin{align*}
	\kappa_{rw}(S)=S^2, \quad \kappa_{ro}(S)=(1-S)^2, \quad \mu_{w}=1,\quad \mu_{o}=5.
\end{align*}
Moreover, the transport equation is
\begin{align}
	S_t+\text{div}(F(S)\bv)=r,\label{transport}
\end{align}
where $F(S)$ is defined as
\begin{align*}
	F(S):=\frac{\kappa_{rw}(S)/\mu_{w}}{\kappa_{rw}(S)/\mu_{w}+\kappa_{ro}(S)/\mu_{o}}
\end{align*}
In this work, $f$ and $r$ are the same, which are defined in \eqref{source}. 

There is an extremely essential quantity called water cut in the above flow problem, which describes the composition of oil-water two-phase flow at the producer. The measurement of water cut can be used to improve oil recovery ratio. Water cut is defined to be $q_w/(q_w+q_o)$, where $q_w:=\int_{\partial D_{\text{out}}}F(S)\bv\cdot \bn dS$ and $q_t:=\int_{\partial D_{\text{out}}}\bv\cdot \bn dS$.

We then describe our scheme to solve the two-phase flow problem. 
To solve \eqref{transport}, we use the finite volume method on the fine grid. In particular, in each fine-scale element $\zeta_i$, we solve the following equation,
\begin{align*}
	|\zeta_i|\frac{S_i^{n+1}-S_i^n}{\Delta t}+\int_{\partial \zeta_i}\hat{S}^n(\bv\cdot \bn)=r_i |\zeta_i|,
\end{align*}
where $|\zeta_i|$ is the area of $\zeta_i$ and $\hat{S}$ is the upwind flux. Besides, $r_i$ is the average value of $r$ on $\zeta_i$. 

To solve the coupling system formed by \eqref{flow_eqn} and \eqref{transport}, we apply a standard Implicit Pressure Explicit Saturation (IMPES) scheme \cite{chen2004improved,chen2006computational}. We define $S^n$ and $\bv^n$ to be the saturation and velocity at the time step $n$. In particular, for a given $S^n$ from the time step $n$, we implicitly solve the flow equation \eqref{flow_eqn} to obtain a velocity solution $\bv^{n+1}$ at the time step $n+1$. Then we explicitly solve $S^{n+1}$ at the time $t_{n+1}$ using \eqref{transport}. In other words, we update the saturation with the velocity solution at the current time step. We give the procedures of the IMPES scheme in Table \ref{tab:impes}. Since a velocity solution computed with mixed finite element method is called a reference velocity solution, we call the corresponding saturation the reference saturation. Similarly, the saturation updated with multiscale velocity bases is denoted by multiscale saturation or approximation of saturation. We will use these definitions in the following discussions.

\begin{table}[htbp!]
	\centering
	\caption{The IMPES scheme}
	\begin{tabular}{c l}
		\hline 
		\hline
		Input:& $S^n$\\
		\hline
		Step 1:& Seek $\bv^{n+1}$ with \eqref{flow_eqn} using $S^n$.\\
		\hline
		Step 2:& Seek $S^{n+1}$ with \eqref{transport} using $S^n$ and $\bv^{n+1}$.\\
		\hline
		Output:& $S^{n+1}$.\\
		\hline
	\end{tabular}
	
	\label{tab:impes}
\end{table}
In a fixed model case, where the permeability field is deterministic, we measure the saturation error by the relative difference in $L^2$ norm between a reference solution and a corresponding approximation; in stochastic case, where different permeability fields are used, we use the relative error of the average saturation computed using MGMsFEM according to the average reference saturation. The mathematical definitions of the concerned errors are given below:
\begin{eqnarray}
	\begin{aligned}
		e_{s,\text{fixed}}:=\frac{\int_{\Omega}|s_{f}-s_{\text{ms}}|^2}{\int_{\Omega}|s_{f}|^2}, \quad
		e_{s,\text{stoc}}:=\frac{\int_{\Omega}|\bar{s}_{f}-\bar{s}_{\text{ms}}|^2}{\int_{\Omega}|\bar{s}_{f}|^2}, \label{e_s}
	\end{aligned}
\end{eqnarray}
where $\bar{s}_{f}$ and $\bar{s}_{\text{ms}}$ are defined to be the average reference saturation and the average multiscale saturation. 

\begin{figure}[htbp!]
	\centering
	\subfigure[$\eta=\frac{1}{4}$]{\includegraphics[width=0.32\textwidth]{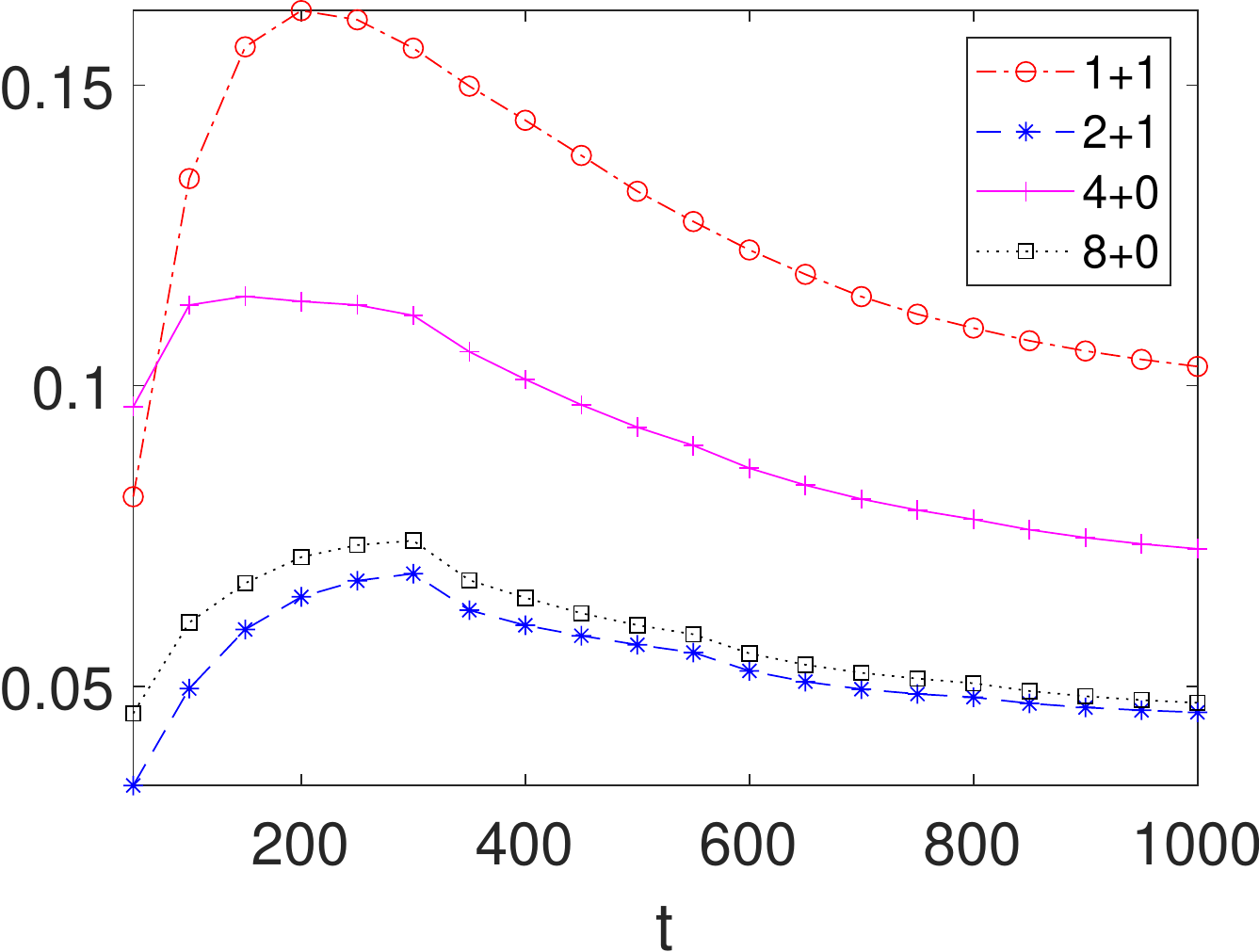}}
	\subfigure[Fixed model $\kappa_4$]{\includegraphics[width=0.32\textwidth]{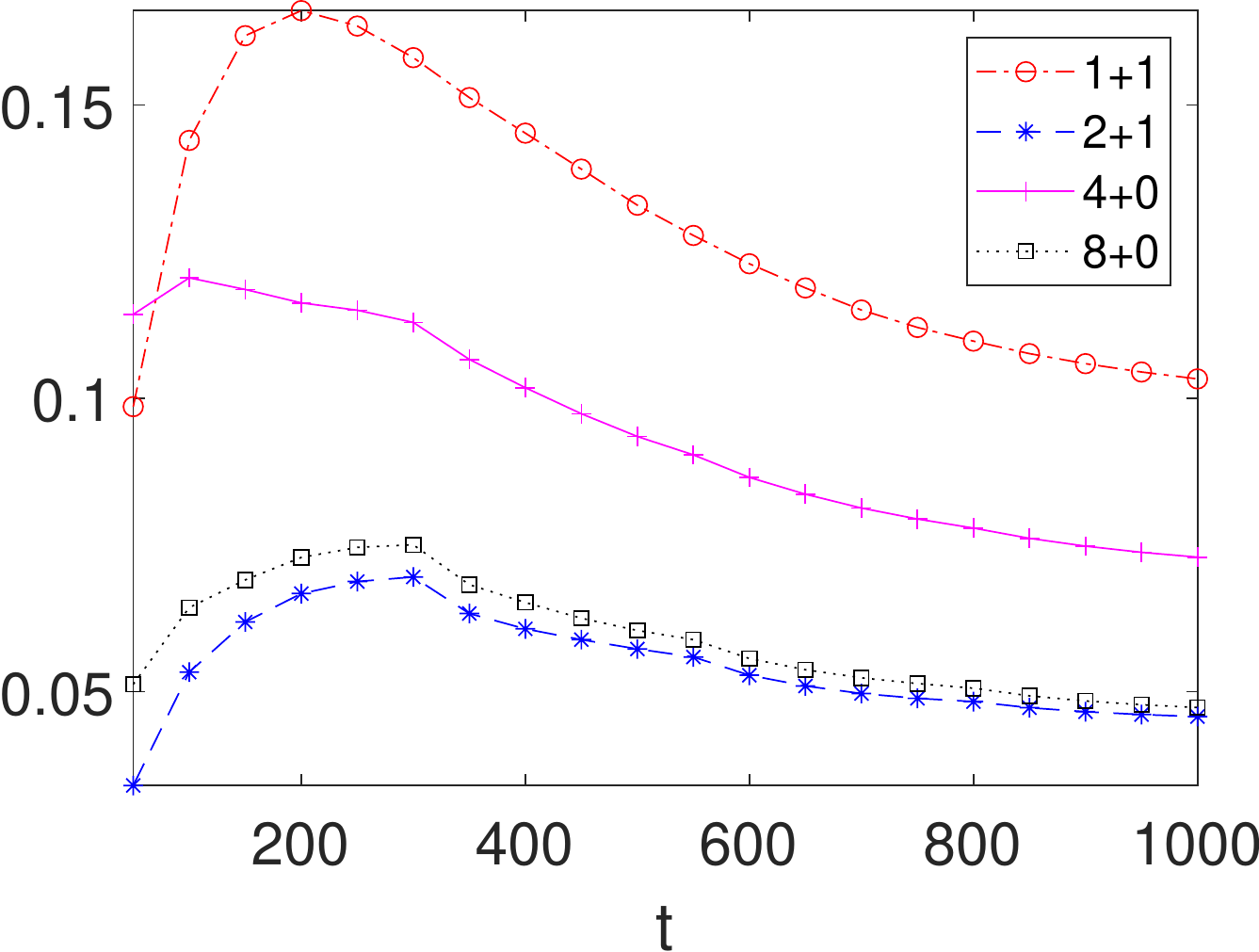}}
	\subfigure[$\eta=\frac{1}{8}$]{\includegraphics[width=0.32\textwidth]{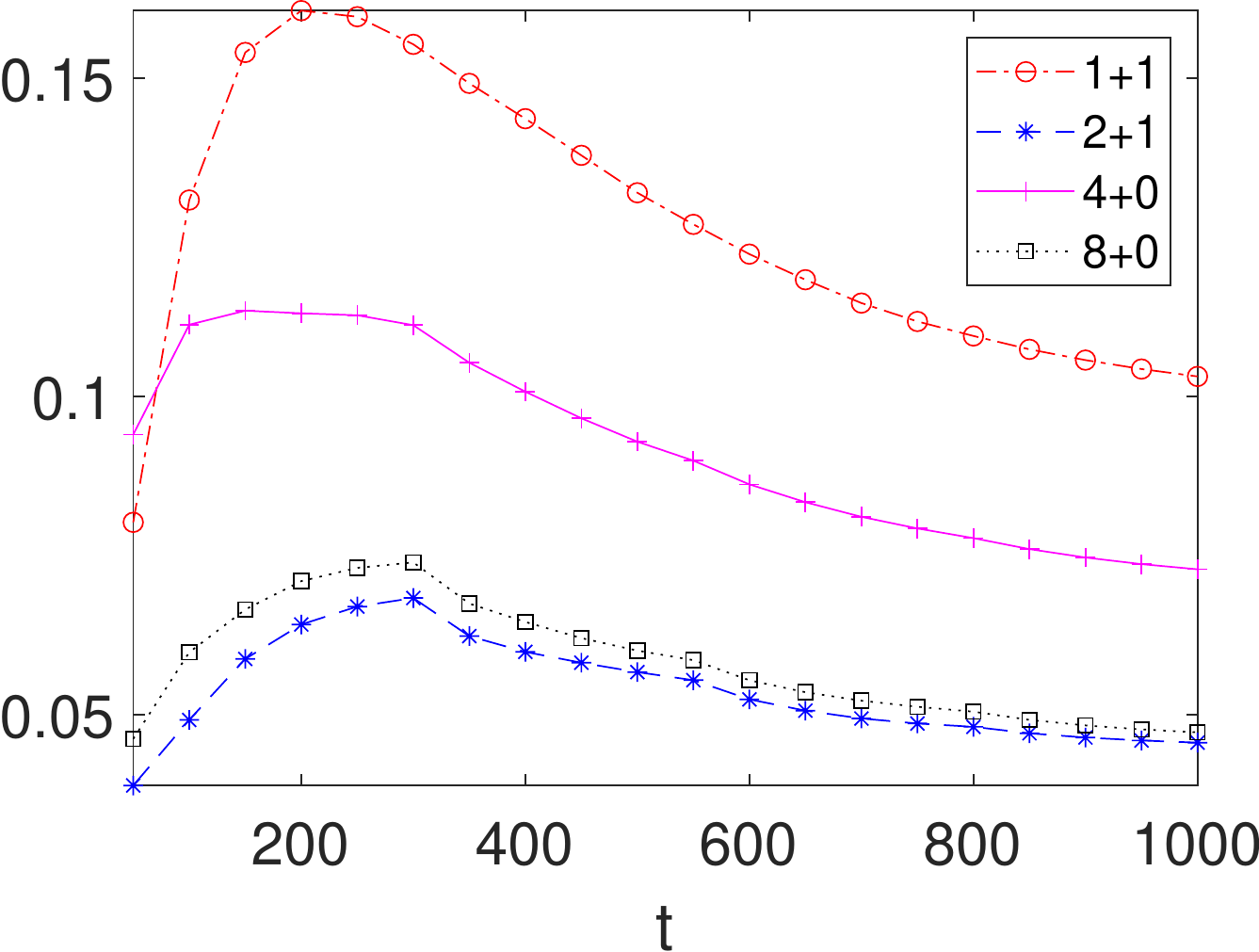}}
	\caption{Comparison of  saturation errors of different numbers of multiscale bases: $``1+1"$, $``2+1"$, $``4+0"$, and $``8+0"$. Left: error of the average saturations (the mean of 100 samples) with $\eta=1/4$;  middle: error of the saturation corresponding to $\kappa_4$; right: error of the average saturations (the mean of 100 samples) with $\eta=1/8$.  $t$ ranges from $0$ to $1000$. }
	\label{secompare}
\end{figure}

In Figure \ref{secompare}, we show the dynamics of saturation errors  associated with different numbers of multiscale bases, i.e. $``1+1"$, $``2+1"$, $``4+0"$ and $``8+0"$. The results are from three groups: a fixed permeability field, two types of random permeability fields with correlation lengths $\eta=\frac{1}{4}$, $\frac{1}{8}$. For each choice of correlation length, we generate 500 sample permeability fields. We then compute the approximation of saturation and the reference saturation for each sample. As shown in \eqref{transport}, one needs to update velocity in a time-marching process. The approximation of velocity is updated by MGMsFEM, while the reference velocity is solved by the mixed finite element method on the fine grid. After we compute 500 approximations and reference solutions of saturation, we compute the average approximation $\bar{s}_{\text{ms}}$ and the average reference saturation $\bar{s}_f$ and compute the $e_{s,\text{stoc}}$ defined in \eqref{e_s}. We place the results from the fixed model in the middle and the two random cases on both sides for a better comparison. Four observations could be obtained. First of all, the performance of multiscale bases are almost the indistinguishable in three cases, which means that the randomness included in the permeability field has a small impact on the approximation effects of multiscale bases. Thus, the extension of MGMsFEM to the concerned two-phase model in random case is shown to be effective.
Secondly, the saturation errors corresponding to different numbers of multiscale bases change in a similar pattern. Specifically, the errors first increase from the initial time to an intermediate time point and then decreases to a level, where it remains steady. The turning points of the four curves are around $t=300$, where we updated the multiscale bases based on the latest saturation. Moreover, residual-driven enrichment is powerful if more than one bases are utilized in each local region in offline stage I. As shown in the graph, the $``1+1"$ case has the largest error among four instances. In particular, the error has exceeded $15\%$ at the peak time. However, once one more offline basis function is added, the error dramatically declined even though only one residual-driven basis is used. More specifically, the $``2+1"$ approximation becomes the most accurate one, which is even better than $``8+0"$ case. Hence, the residual-driven bases are powerful when sufficient offline bases are used. Last but not least, the effect of using more offline bases is still apparent when no residual-driven bases are involved. In Figure \ref{secompare}, an evident improvement can be noticed from $``4+0"$ to $``8+0"$. More specifically, the error is almost halved when the number of bases is doubled. Consequently, even though the residual-driven bases are powerful, one can not disregard the significance of offline bases.

\begin{figure}[htbp!]
	\centering
	\subfigure[fine saturation at $t=250$ ]{\includegraphics[width=0.32\textwidth]{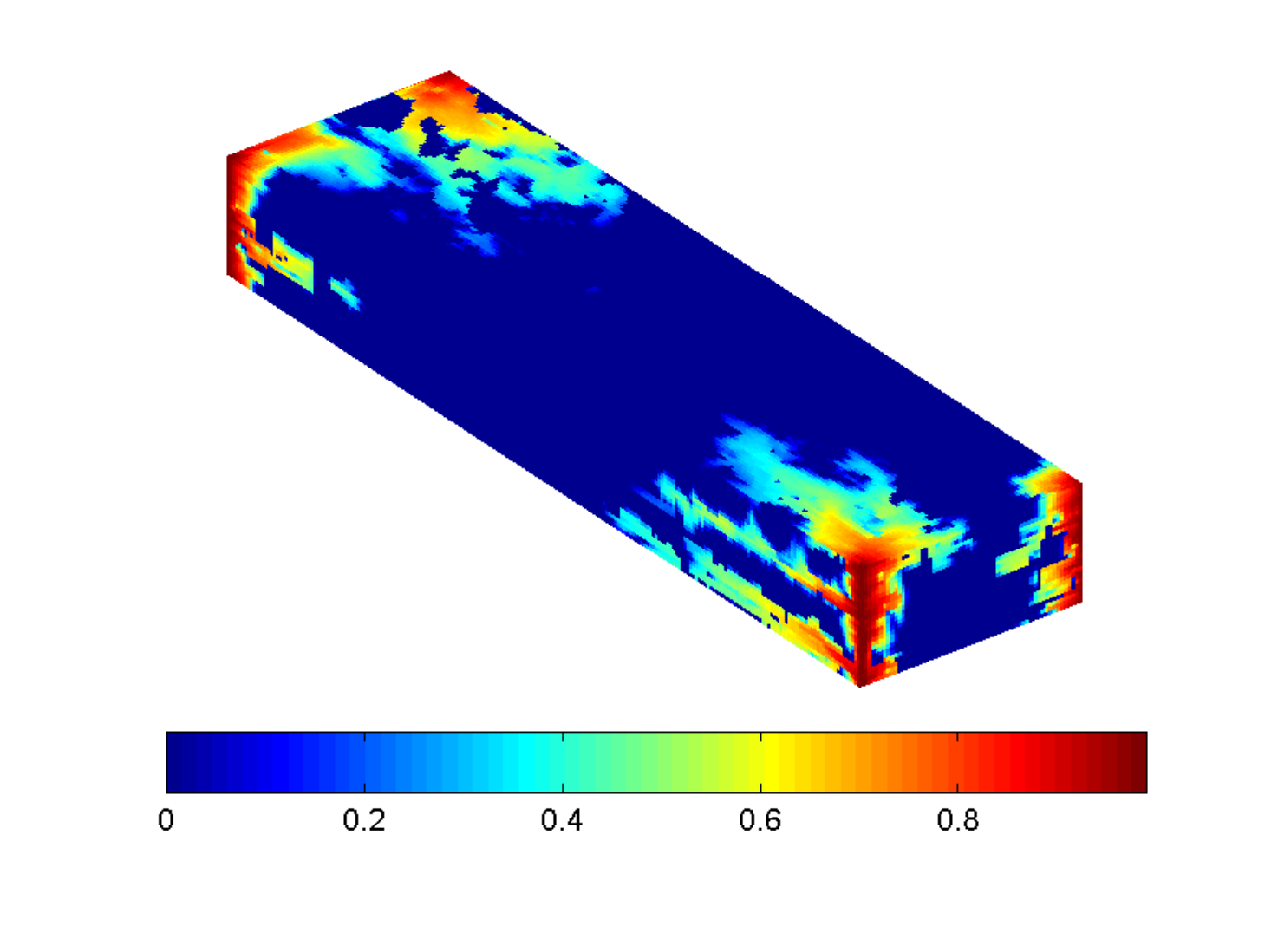}}
	\subfigure[ms saturation with 2+1 bases at $t=250$ ]{\includegraphics[width=0.32\textwidth]{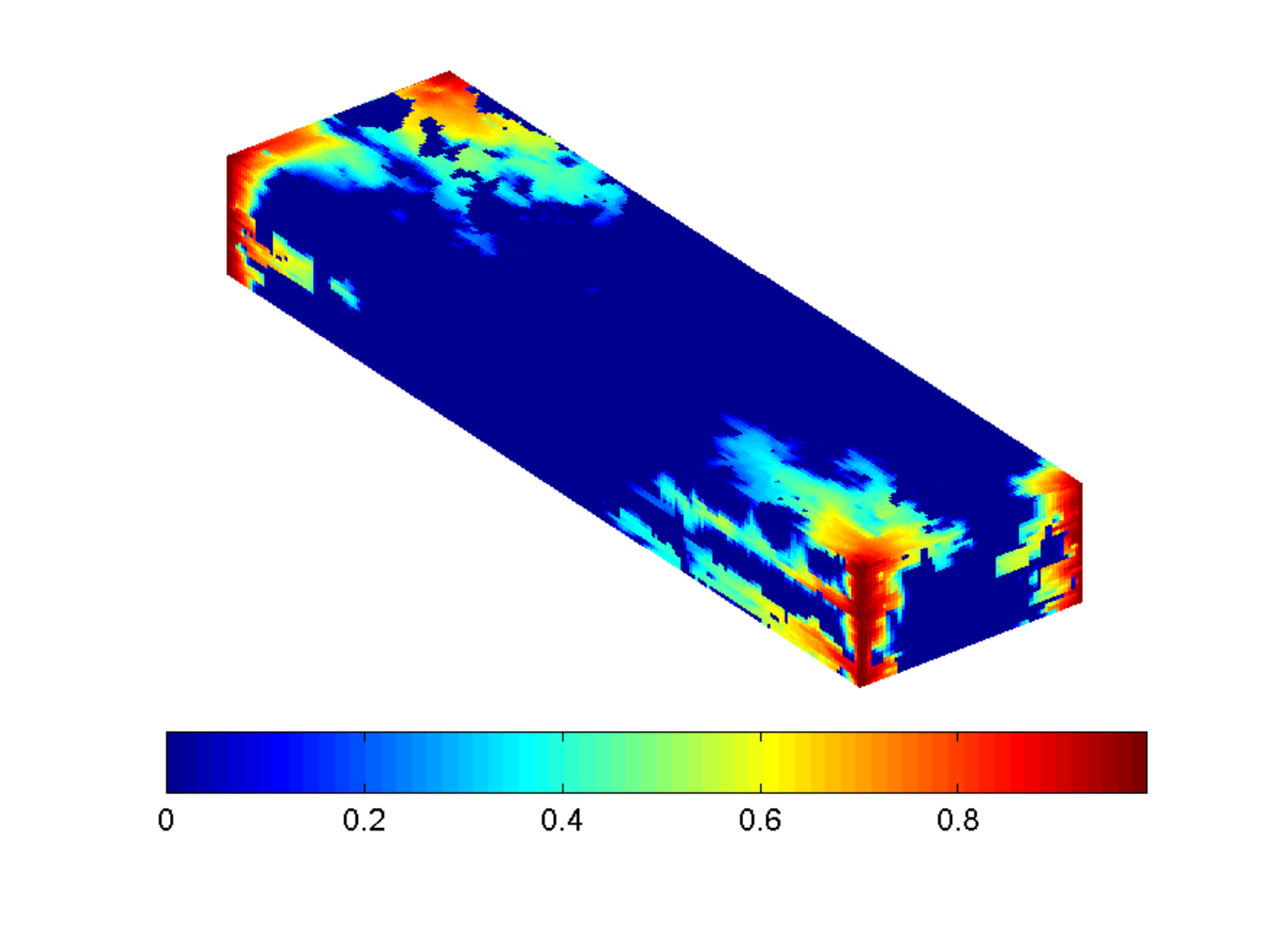}}	
	\subfigure[absolute difference at $t=250$ ]{\includegraphics[width=0.32\textwidth]{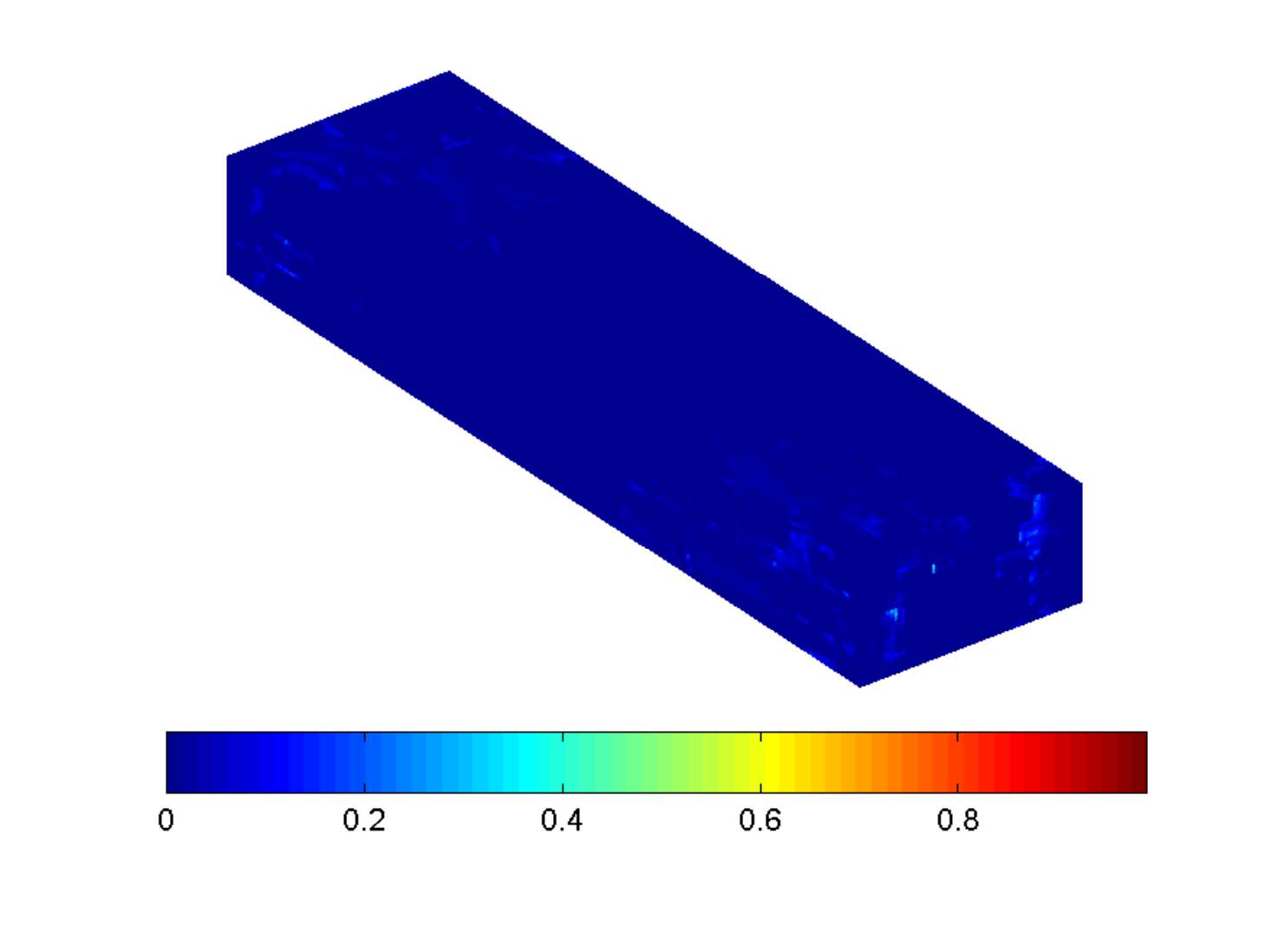}\label{e_t250}}	
	\subfigure[fine saturation at $t=750$ ]{\includegraphics[width=0.32\textwidth]{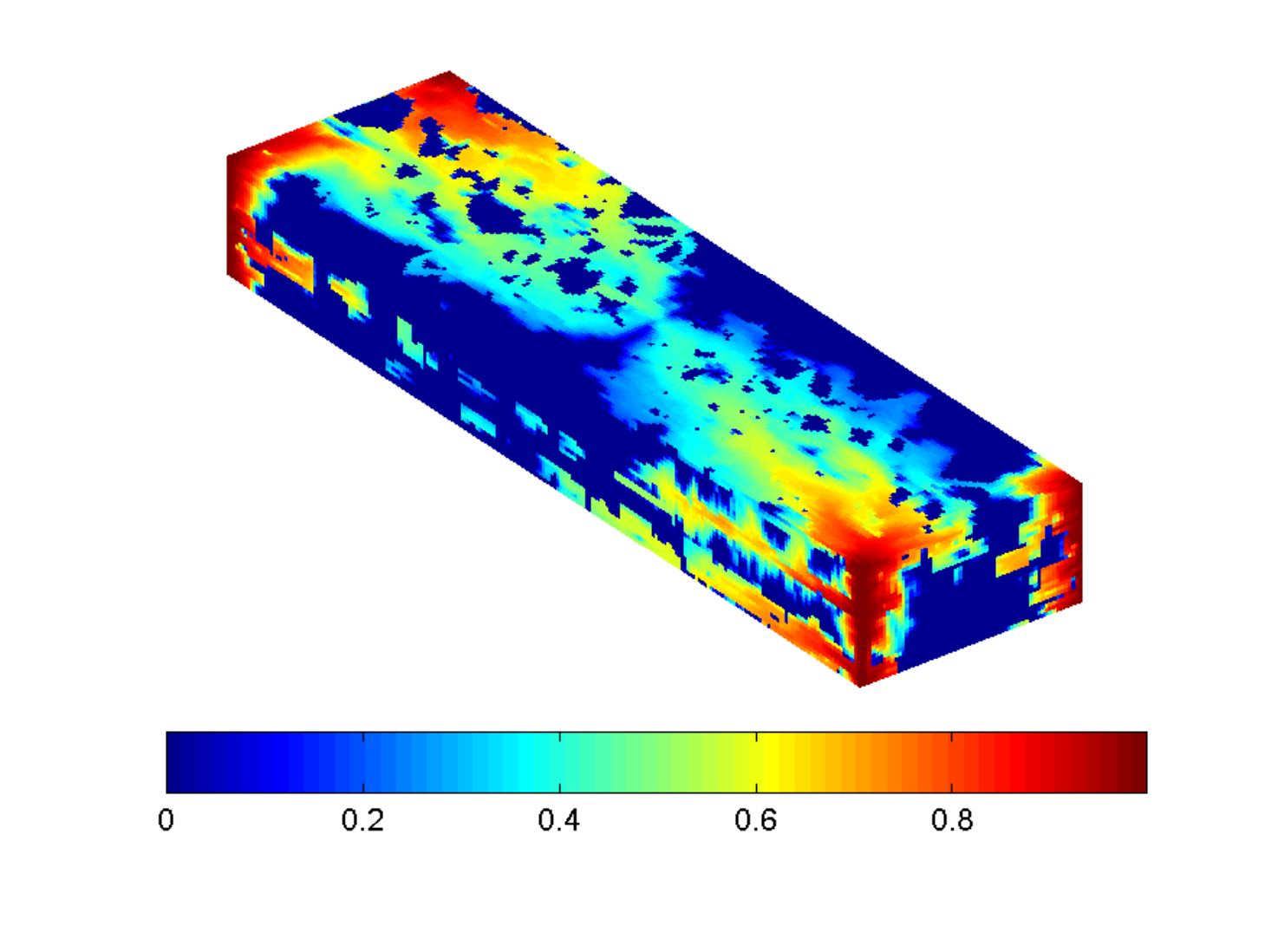}}
	\subfigure[ms saturation with 2+1 bases at $t=750$ ]{\includegraphics[width=0.32\textwidth]{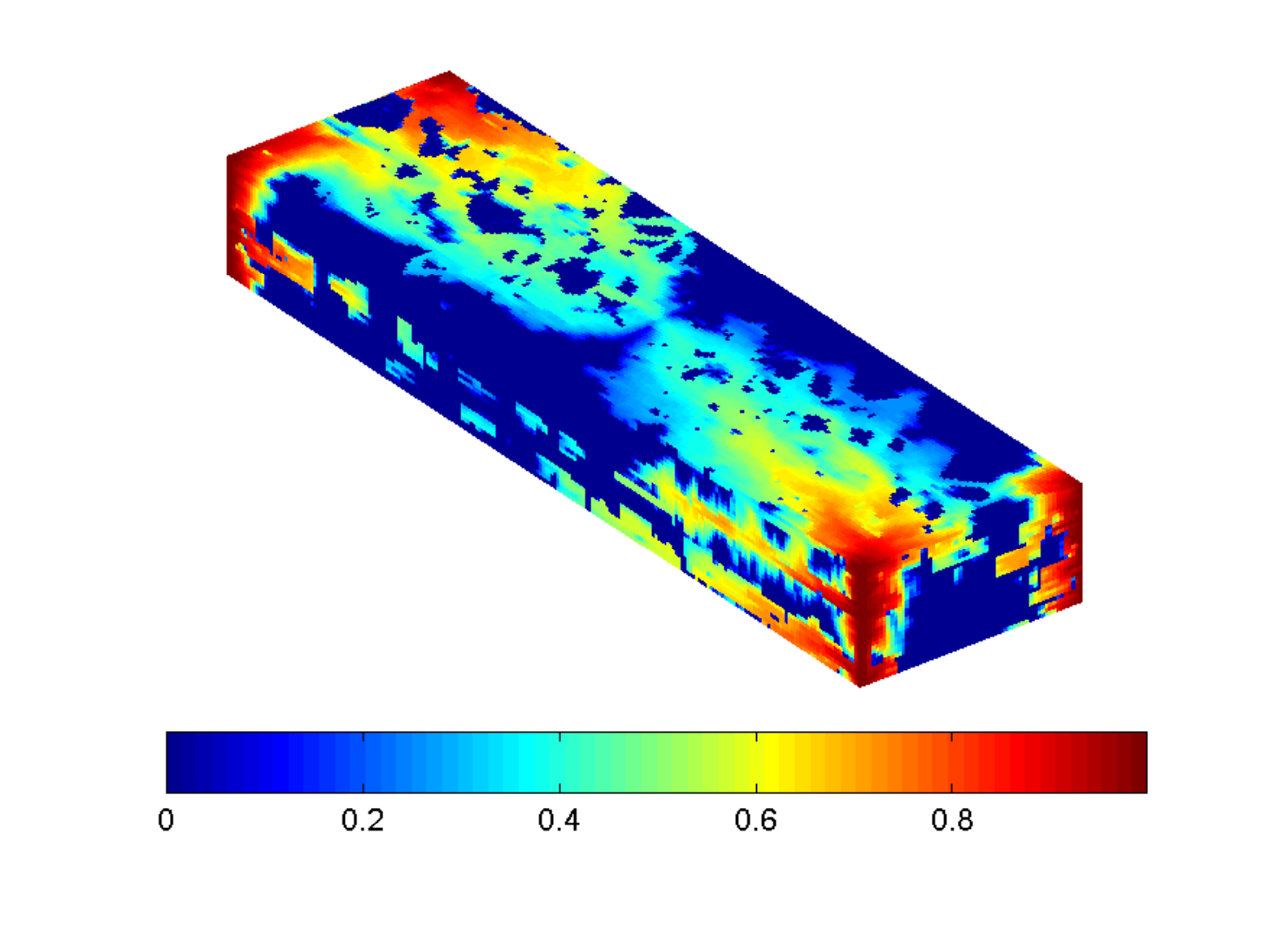}}	
	\subfigure[absolute difference at $t=750$ ]{\includegraphics[width=0.32\textwidth]{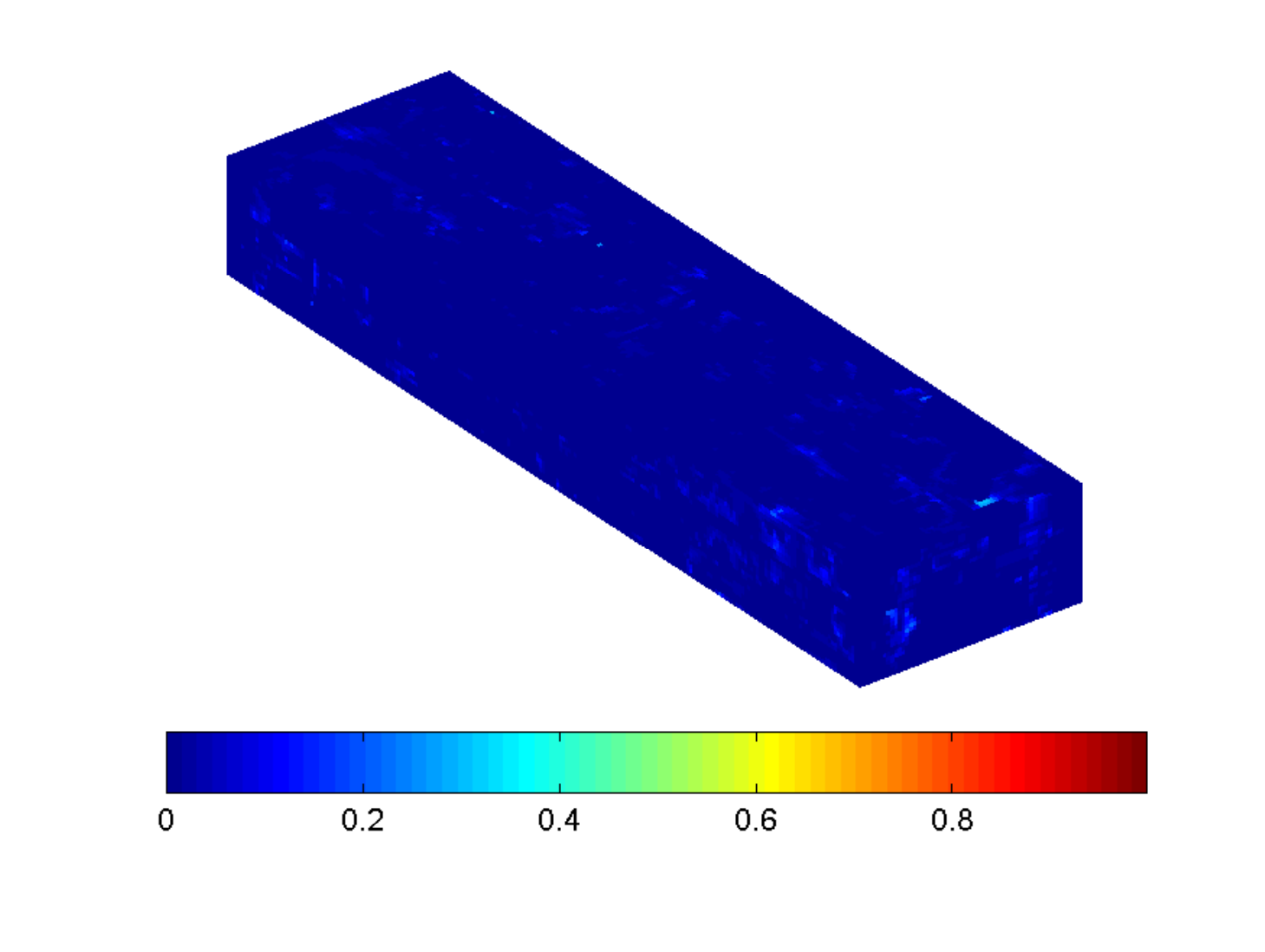}\label{e_t750}}	
	\caption{Comparison of the  average fine saturation and the average multiscale saturation with $``2+1"$ bases at $t=250$ and $t=750$. 100 samples are used, where the mean permeability field is $\kappa_4$.}
	\label{fig:sat_2_1}
\end{figure}

In Figure \ref{fig:sat_2_1}, we compare an average reference saturation with an average  multiscale saturation associated with $``2+1"$ bases respectively, where the averages are computed among 100 samples. To better visualize the comparison, the absolute difference between the above two concerned quantities is shown in the last column. We demonstrate two time steps: $t=250$ and $t=750$.  The reference saturation is shown on the left. The flows move from the four vertical edges and extend to the center. One can see that the simulation is faster where the corresponding permeability shows a greater variation. In particular, the flows starting from the lower-left and upper-right corners shift further than the other two corners within the same time. The simulations of $t=250$ and $t=750$ are displayed in the first and second rows. The flows have not  intersected at the center until $t=750$. At $t=250$, the multiscale saturation is a relatively accurate approximation of the reference saturation because most detailed information is shown to be preserved in the multiscale saturation. In other words, it is difficult to find noticeable differences between them. At $t=750$, more details are displayed. Overall, the multiscale saturation is indistinguishable from the reference. However, there is an observable distinction at the lower-right corner, where the permeability shows strong discontinuity in the subfigure \ref{spe10last30}. The permeability shown in Figure \ref{fig:models} is highly discontinuous at this region, which reveals some connections between the permeability field and the saturation. In other places, there are barely big deviations between the approximation and the fine saturation, which shows the effect of multiscale bases. 

\begin{figure}[htbp!]
	\centering
	\subfigure[$\eta=\frac{1}{4}$ ]{\includegraphics[width=0.48\textwidth]{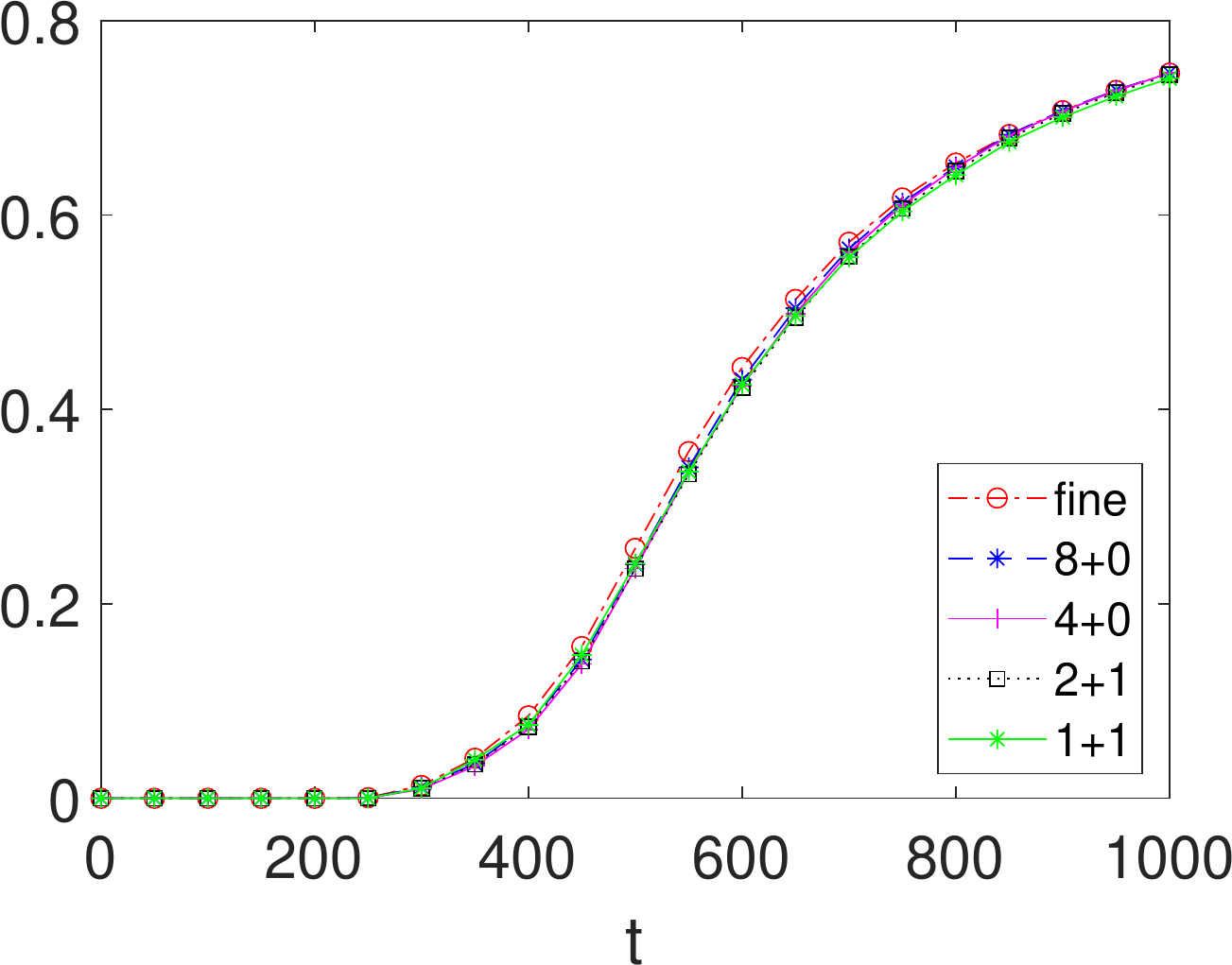}}
	\subfigure[$\eta=\frac{1}{8}$]{\includegraphics[width=0.48\textwidth]{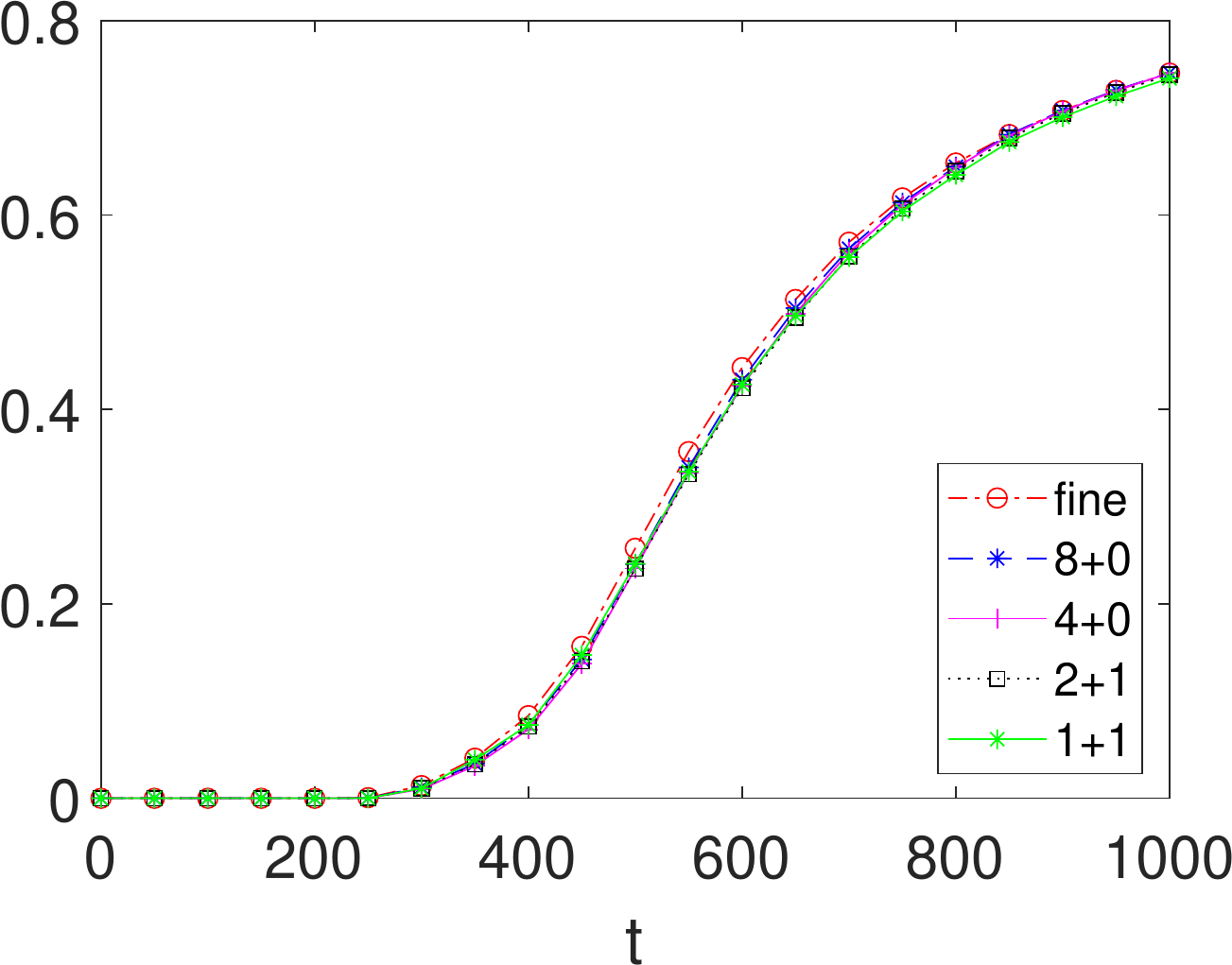}}
	\caption{Average water-cut curves computed with different methods under two choices of correlation lengths, $\eta=1/4,1/8$. Five curves are corresponding to MFEM and MGMsFEM with four numbers of multiscale bases (``8+0'', ``4+0'', ``2+1'' and ``1+1''). 100 samples are used, where the mean permeability field is $\kappa_4$. }	
	\label{fig:watercut}
\end{figure}

In Figure \ref{fig:watercut}, we show water-cut curves corresponding to a reference and four approximations. The reference water cut is computed by MFEM on a fine grid hence we use ``fine'' to denote this curve. Four approximation water cut curves are solved with MGMsFEM, where the numbers of multiscale bases are  $``1+1"$, $``2+1"$, $``4+0"$, and $``8+0"$. We demonstrate the comparisons with $\eta=\frac{1}{4}, \frac{1}{8}$ from $t=0$ to $1000$. From the two figures, one can observe that all the water cuts first increase slowly then fast and tend steady, which shows the dynamics of fraction of water at the producer. The differences between the approximations and the references are indistinguishable, which holds for all the approximations. Water cut is indeed a function of saturation at specific points, which is the reason that the approximations with only bases in offline stage I are nearly as accurate as the those computed with bases in both two offline stages.  It is worth mentioning that even though the approximation bias can be noticed  at the middle time, the errors tend to vanish as the water tends to be saturated.

\section{Conclusion}
In this work, our objective is to provide a computationally efficient method for stochastic groundwater flow problems. In particular, we apply a mixed generalized multiscale finite element method to a single-phase flow equation and a two-phase flow equation coupled with a transport equation, where randomness is included in permeability fields of the above two problems. With a well-chosen permeability field, one can construct a set of locally-defined multiscale bases during two stages in training, which will be further used to solve equations in a coarse mesh with an arbitrary permeability field in a testing stage. 

We utilize a mixed formulation of the flow equation to achieve a conservation of mass in local regions, which is essential in flow problems. The multiscale space in offline stage I is constructed based on a much larger snapshot space, where snapshot bases are resulted from different boundary conditions. Using a set of well-designed spectral problems, we can distill some important modes from the snapshot space. Since the constructions are performed independently in each local patch, one can use a parallel computing to achieve higher efficiency. Based on a set of local residuals derived from the equation, some residual-driven bases are generated to enhance the accuracy of the solutions to the flow equation and further decrease the errors in saturation solutions to the transport equation in the two-phase model. We emphasize that although the residual-driven construction is conducted in a fine-grid mesh, the computations are restricted in some local patches. In other words, the computation cost exists in the residual-driven enrichment is limited and the constructions of residual-driven bases can also be implemented in parallel. Overall, the final multiscale space is spanned by multiscale bases in two offline stages. Because the total number of bases is remarkably reduced compared with the previous snapshot space, the computation is significantly accelerated. Hence, our method can be employed to speed up the uncertainly quantification problems in subsurface flows, where a large number of samples are investigated. What's more, the proposed method does not require any interpolation in the spatial domain and can be applied to solve the equation associated with any realizations sampled from a large stochastic space.

We demonstrate numerical simulations for a single-phase flow and a two-phase flow. More specifically, we use stochastic high-contrast models and SPE models. For each model, both two-dimensional and three-dimensional cases are investigated. For the single-phase flow, our method focuses on reducing the velocity errors, while for the latter one, a good accuracy in the saturation solutions is the target. Our results show that the multiscale space generated with a training field can be applied to different sample fields as well as different sources. Furthermore, some residual-driven bases can improve the overall accuracy. In particular, in the single-phase flow, we show test errors corresponding to sample fields and a source that are different from what are used in the training stage. The errors are small if sufficient bases are used. The use of residual-driven bases can remarkably increase the accuracy and efficiency especially in 3D case.
For the two-phase flow model, we show a comparison of reference saturations and approximations of saturation at two specific time points. The  approximations can well capture most  details in the references, which means the differences between these two saturations are almost difficult to distinguish. Besides, we also compute water cut in the context of two-phase flow. Since water cut is a  deterministic function of saturation, the approximations of water cut based on multiscale space are pretty accurate compared to the reference water cut, which is similar to approximation reference saturations. Besides, we would like to remark that the proposed method can be extended to a more general case, where different types of grids (not necessarily Cartesian grids) are used for a spatial discretization.

Generally, the proposed method is promising in dealing with groundwater flow problems, which is guaranteed by theoretical proofs and verified in numerical results. In the future, there are still some interesting and non-trivial problems to explore. For example, we may consider a more sophisticated random space and combining our method with a stochastic collocation method \cite{babuvska2007stochastic}. Besides, we can further apply the proposed algorithm to some inverse problems, where the associated forward process can solved by our method. Last but not least, this method can serve as a base to some non-intrusive methods like deep learning. One can combine this method with deep learning as our recent work \cite{Wang2021DL}.

\section*{Acknowledgement}

The research of Eric Chung is partially supported by the Hong Kong RGC General Research Fund (Project numbers 14304719 and 14302620) and CUHK Faculty of Science Direct Grant 2020-21.

\bibliographystyle{plain}
\bibliography{references}
\end{document}